\documentclass[pdflatex,sn-mathphys-num]{sn-jnl}

\usepackage{indentfirst}
\usepackage{graphicx}%
\usepackage{multirow}%
\usepackage{amsmath,amssymb,amsfonts}%
\usepackage{amsthm}%
\usepackage{mathrsfs, dsfont, bbm, bm} 
\usepackage{epstopdf}
\usepackage{epsfig}
\usepackage{caption}
 \usepackage{subfigure}
\usepackage{exscale, relsize}  
\usepackage{chngpage}
\usepackage{mathrsfs}
\usepackage{multicol}
\usepackage{dsfont}
\usepackage{exscale}
\usepackage{tikz-cd}
\usepackage{relsize}
\usepackage{amssymb}
\usepackage{hyperref}
\usepackage{url}
\usepackage{tikz-cd}
\usepackage{comment}
\usepackage{setspace}
\usepackage[misc,geometry]{ifsym} 
\usepackage{color, xcolor}%
\usepackage{textcomp}%
\usepackage{manyfoot}%
\usepackage{booktabs}%
\usepackage{algorithm}%
\usepackage{algorithmic}%
\usepackage{listings}%

\theoremstyle{thmstyleone}%
\newtheorem{theorem}{Theorem}[section]
\newtheorem{definition}{Definition}[section]
\newtheorem{example}{Example}[section]

\newtheorem{remark}{Remark}[section]
\newtheorem{lemma}{Lemma}[section]

\numberwithin{equation}{section}%
\numberwithin{table}{section}%
\numberwithin{figure}{section}

\def\3bar{{|\hspace{-.02in}|\hspace{-.02in}|}}

\renewcommand\div{\operatorname{div}}

\def\d{\text{d}}

\graphicspath{{../figures/}}
\begin{document}
	\title[A enriched Galerkin method for the Navier--Stokes equations]{A pressure-robust and  parameter-free enriched Galerkin method for the Navier--Stokes equations of rotational form}
    
\author[1]{\fnm{Shuai} \sur{Su}}\email{shuaisu@bjut.edu.cn}

\author[1]{\fnm{Xiurong} \sur{ Yan}}\email{yanxr@emails.bjut.edu.cn}

\author*[2]{\fnm{Qian} \sur{ Zhang}}\email{qzhang25@jlu.edu.cn}

\affil[1]{\orgdiv{School of Mathematics, Statistics and Mechanics}, 
\orgname{Beijing University of Technology}, 
\orgaddress{\city{Beijing}, \postcode{100124}, \country{China}}}

\affil*[2]{\orgdiv{School of Mathematics}, 
\orgname{Jilin University}, 
\orgaddress{\city{Changchun}, \postcode{130012}, \country{China}}}

	

\abstract{In this paper, we develop a novel enriched Galerkin (EG) method for the steady incompressible Navier–Stokes equations in rotational form, which is both pressure-robust and parameter-free. The EG space employed here, originally proposed in \cite{su2025novel}, differs from traditional EG methods: it enriches the first-order continuous Galerkin (CG) space with piecewise constants along edges in two dimensions or on faces in three dimensions, rather than with elementwise polynomials. Within this framework, the gradient and divergence are modified to incorporate the edge/face enrichment, while the curl remains applied only to the CG component, an inherent feature that makes the space particularly suitable for the rotational form. The proposed EG method achieves  pressure robustness through a velocity reconstruction operator. We establish existence, uniqueness under a small-data assumption, and convergence of the method, and confirm its effectiveness by numerical experiments.
}	
	






\maketitle
\section{Introduction}\label{sec1}
The Navier--Stokes (NS) equations are a fundamental model for incompressible Newtonian fluid dynamics. They play a crucial role in various applications, such as fluid flow in pipes and channels, aerodynamic flows around airplane wings, weather and climate studies, and circulation of blood, to name a few. 
The nonlinearity in the NS equations can be written in several forms, including the convective form, the skew symmetric form, and the rotational form. While these forms are equivalent at the continuous level, they lead to discretizations with varying algorithmic costs, conserved quantities, and levels of approximation accuracy. 
Although the convective form is a common choice in the finite element discretization of the NS equations, the rotational form provides several advantages that make it a compelling alternative. In particular, it better preserves key physical properties such as helicity and enstrophy conservation, often leads to improved stability in numerical simulations, facilitates the development of efficient iterative solvers, and, compared to the skew-symmetric form, is typically less expensive to compute -- making it a natural starting point for our investigation \cite{layton2009accuracy, rebholz2007energy, lube2002stable, olshanskii2002low,olshanskii1999iterative}. 
In this paper, we consider the steady-state NS equations  of rotational form on a  bounded and connected domain  $\Omega$ in $\mathbb{R}^d (d=2,3)$: 
\begin{subequations}\label{eq:n1}
\begin{align}
-\nu \Delta {\bm u}+(\nabla\times\bm u) \times {\bm u}+\nabla p & ={\bm f} \quad \text { in } \Omega, \label{eq:n1a}\\
\nabla \cdot {\bm u} & =0 \quad \text { in } \Omega, \label{eq:n1b}\\
{\bm u} & =0 \quad \text { on } \partial \Omega, \label{eq:n1c}
\end{align}
\end{subequations}
 where  $\bm {u}$ is  the velocity,   $p$ is the pressure, 
 ${\bm f}$ is a given body force,  and $\nu>0$ is the viscosity of the fluid. 
Here the pressure $p$  is called the Bernoulli pressure satisfying the equation $p=p^{\mathrm{kin}}+\frac{1}{2}|\bm u|^2$ with the kinematic pressure $p^{\mathrm{kin}}$. When $d=2$, $(\nabla\times\bm u)\times \bm u =(\partial_xu_2-\partial_yu_1)\left(-u_2,u_1\right)^{\mathrm T}$. 

In the community of mixed finite element methods for the (Navier--)Stokes equations, to ensure stability, it is essential to use velocity and pressure pairs satisfying the inf-sup condition. A variety of finite element pairs have been constructed in previous works to meet this requirement, which include classical finite element \cite{boffi2013mixed,Girault2012Finite} as well as variations like discontinuous Galerkin (DG) \cite{cockburn2009equal,cockburn2005locally}, hybrid discontinuous Galerkin (HDG) \cite{cesmelioglu2017analysis}, hybrid high-order (HHO) method \cite{botti2019hybrid,di2018hybrid}, and weak Galerkin (WG) \cite{hu2019weak,zhang2018weak}, etc.
Nevertheless, even in some inf-sup stable mixed finite element methods, such as the nonconforming Crouzeix-Raviart element and the conforming Taylor-Hood element, a particular type of nonrobustness is evident \cite{john2017divergence}. This nonrobustness is characterized by the velocity error depending on a pressure error contribution, expressed as \(\frac{1}{\nu}\inf_{q_h\in W_h}\|p-q_h\|\) with the discrete pressure space \(W_h\) and viscosity coefficient $\nu$.
This issue is primarily related to how the mass conservation is maintained (how the divergence constraint is discretized), rather than stemming from the nonlinearity or dominant convection of the equations. Consequently, this lack of robustness is often described as poor mass conservation, and is traditionally mitigated by grad–div stabilization \cite{olshanskii2004grad,case2011connection,olshanskii2009grad}. The desired robustness is termed pressure-robustness, implying that the velocity error remains unaffected by the pressure error.

The NS equations of the rotational form  require the use of the Bernoulli pressure, which is particularly challenging to solve at high Reynolds numbers, and hence leads to inferior velocity discretization in schemes lacking pressure-robustness \cite{layton2009accuracy}. Therefore, it is critical to develop pressure-robust numerical methods for the NS equations of this form. Recently, a novel technique of using a velocity reconstruction operator has been proposed \cite{linke2014role}, which maps the velocity test functions to an $H(\div)$-conforming finite element space.
This technique has shown promise in achieving pressure-robustness in certain mixed finite element methods for the rotational form \cite{linke2016pressure,mu2023pressure,Daniel2020,yang2022analysis}.
 Specifically, Linke et al. \cite{linke2016pressure} and Yang et al. 
 \cite{yang2022analysis} have modified the existing inf-sup stable conforming finite element methods by incorporating this velocity reconstruction technique to ensure pressure-robustness. Notably, Yang et al. employed an equivalent formulation under $H^1$ regularity rather than the rotational form directly. Similarly, the WG method of Mu \cite{mu2023pressure} and the HHO method of Quiroz et al. \cite{Daniel2020} are also based on such an equivalent formulation. Moreover, both WG and HHO discretizations generally require higher computational costs than conforming finite element methods for the same level of accuracy. The aim of this work is to develop a pressure-robust and computationally efficient finite element method that directly targets the rotational form of the NS equations.

Recently, the enriched Galerkin (EG) method has attracted considerable attention due to its high efficiency and ease of implementation. It was first introduced in \cite{SL2009} to solve a second-order elliptic problem and was shown to be locally mass conservative.
The basic idea of the EG method
is to enrich the continuous Galerkin (CG) finite element space with a DG space and use it in the DG
formulation. This approach can maintain the desirable features of the DG method at a computational
cost comparable to that of the CG method. To date, the EG method  
has been successfully applied in various problems, such as elliptic and parabolic problems in porous media 
\cite{lee2016locally}, two-phase flow 
\cite{lee2018enriched}, the shallow-water equations \cite{hauck2020enriched}, the Stokes problem \cite{chaabane2018stable,yi2022enriched,hu2024pressure,lee2024low}, and linear elasticity \cite{yi2022locking,peng2024locking,su2024parameter}. 
In particular, a pressure-robust EG method has been proposed in \cite{hu2024pressure} for the Stokes problem by applying the velocity reconstruction technique. However, it remains uncertain whether this method can be effectively adapted to solve the NS equations of rotational form.
In our previous work \cite{su2025novel}, we (Su and Q. Zhang, together with collaborators Tong and M. Zhang) introduced a novel EG space for solving linear elasticity problems in both two and three dimensions.
Unlike traditional EG methods that enrich the CG space with a DG space defined on elements, our approach enriches the first-order CG space with a DG space of piecewise constants along edges in two dimensions (2D) or on faces in three dimensions (3D). 
This DG enrichment acts as a correction to the normal component of the CG space, facilitating the establishment of the inf-sup condition compared to the classical EG space in \cite{yi2022locking}.
By incorporating these edge or face corrections, we define both a  modified divergence and a modified gradient. However, the curl operator remains applied solely to the CG component, as the DG enrichment involves only the normal component. 
This distinctive feature makes the proposed EG space particularly well-suited for the rotational form of the NS equations, where the nonlinearity involves the curl of the velocity field.
Motivated by these advantages, we extend the new EG space in \cite{su2025novel} to the NS equations \eqref{eq:n1}. 
The EG space incorporates edge or face components, sharing a similarity with spaces used in the HDG, HHO, and WG methods. This similarity inspires us to develop our EG scheme within the framework of these method, yielding a naturally parameter-free scheme. 
To achieve pressure-robustness, we utilize the velocity reconstruction operator \(\mathcal{R}\) \cite{linke2014role} on both the right-hand side and the nonlinear term of our model. Specifically, we discretize the right-hand side \((\bm{f}, \bm{v})\) and the nonlinear term \(((\nabla\times \bm{u} )\times \bm{u}, \bm{v})\) in the variational formulation by \((\bm{f}, \mathcal{R} \bm{v}_h)\) and \(((\nabla\times \bm{u}_h) \times \mathcal{R} \bm{u}_h, \mathcal{R} \bm{v}_h)\), respectively. Theoretically, we establish the well-posedness of the newly proposed EG method and provide rigorous pressure-robust error estimates. 
It is worth noting that our approximation space is essentially equivalent to that of \cite{li2022low,li2024analysis}, which enriches the first-order CG space with the lowest-order Raviart–Thomas elements to solve the Stokes and linear elasticity problems, and is further applied to the NS equations in \cite{ahmed2024inf}. However, the resulting formulation is fundamentally different: our method is developed within a parameter-free framework and is directly based on the rotational form of the NS equations.

The remainder of this paper is organized as follows. In Section \ref{sec3}, we construct our new EG space and present our PR\&PF-EG method for the NS equations of rotational form. In Section \ref{sec:Theor}, we establish a rigorous theoretical analysis of the newly proposed EG method. We prove the existence of solutions, and demonstrate uniqueness under a small data condition, and derive pressure-robust error estimates. Several numerical experiments are conducted in Section \ref{sec5} to  validate the theoretical results. Finally, some concluding remarks are given in section \ref{conclusion}.

Throughout this paper, we denote by $H^s(D)$  the Sobolev space  with the norm $\|\cdot\|_{s,D}$ for a bounded
Lipschitz domain  $D\subset \mathbb R^n$, $n=1,2,3$ and a real number $s\geq 0$. The space $H^0(D)$ coincides with $L^2(D)$ with the $L^2$-inner product denoted by $(\cdot,\cdot)_D$. When \(D = \Omega\) or \(s = 0\), we omit the subscript \(D\) or \(s\), respectively, whenever there is no confusion. We denote by $P_{\ell}(D)$ the space of $\ell$-th order polynomials on $D$. 
These notations are generalized to vector- and tensor-valued Sobolev spaces.
We use $C$ to denote a generic and positive constant, which is independent of the mesh size $h$.
\section{A pressure-robust and parameter-free enriched Galerkin method}\label{sec3} 
In this section, we  present our PR\&PF-EG   method for the
NS equations. To this end,
we use \(\mathcal{T}_h\) to denote a shape-regular partition of the domain \(\Omega\), consisting of triangles in two dimensions or tetrahedrons in three dimensions. Denote by \(\mathcal{E}_h\) the set of all edges/faces in \(\mathcal{T}_h\).
For  a generic element $T \in \mathcal{T}_h$,  let $h_T$ be diameter  of $T$, and
$h=\max _{T \in \mathcal{T}_h} h_T$ be the mesh size of $\mathcal{T}_h$. 
We denote by $P_{\ell}(\mathcal T_h)$ the space of piecewise $\ell$-th order polynomials over $\mathcal T_h$. 

For completeness, we provide an introduction to the EG space in \cite{su2025novel}.
We first consider the vector-valued linear CG finite element space
\begin{equation*}
 \mathrm{CG}=\left\{\boldsymbol{v} \in [H^1(\Omega)]^d:\left.\boldsymbol{v}\right|_T \in\left[P_1(T)\right]^d \hbox{ for all }  T \in \mathcal{T}_h\right\} \\
\end{equation*}
and the DG space consisting of   piecewise constant functions over $\mathcal E_h$
\begin{equation*}
\mathrm{DG}=\left\{v \in L^2\left(\mathcal{E}_h\right):\left.v\right|_e \in P_0(e) \hbox{ for all }   e \in \mathcal{E}_h\right\}. 
\end{equation*}
Then  the EG space $\bm V_h$ is defined by
$$
\bm V_h :=\left\{\boldsymbol{v}_h=\left\{\boldsymbol{v}_0, v_b\right\}: \boldsymbol{v}_0 \in \mathrm{CG} \text { and } v_b \in \mathrm{DG}\right\} .
$$
Here we remark that the component $v_b|_e$ serves as a correction to $\frac{1}{|e|} \int_e \boldsymbol{v}_0 \cdot \boldsymbol{n}_e \mathrm{~d} s$, where $\boldsymbol{n}_e$ represents the assigned unit normal vector to the edge $e$. 
Our EG space shares certain characteristics with the spaces
used in the HDG, HHO, and WG methods. The key distinction, however, lies in the usage of the
CG space for $\boldsymbol{u}_0$, as opposed to the discontinuous piecewise polynomial space employed in these methods.

For the pressure $p$, we simply use the piecewise constant function space
$$
W_h=\left\{ q \in L^2(\Omega): \left.q\right|_T \in P_{0}(T)\text{ for all } T \in \mathcal{T}_h\right\} .
$$

On one hand,  to avoid tuning penalty parameter, we define a modified gradient and a modified divergence for any $\boldsymbol{v} =\{\boldsymbol{v}_0,v_b \}\in \bm V_h $ by replacing the normal component in the integration by parts formula with the enriched component $v_b$.

\begin{definition}[modified gradient \cite{su2025novel}]
The modified gradient operator for $\boldsymbol{v} \in \bm V_h $ is defined as  $\nabla_{m} \boldsymbol{v} \in\left[P_0(\mathcal T_h)\right]^{d \times d}$ such that for all $\sigma \in\left[P_0(T)\right]^{d \times d}$ and  $T\in\mathcal T_h$,
\begin{equation}\label{eq:grad}
\begin{aligned}
\left(\nabla_{m} \boldsymbol{v},\sigma\right)_T=\left\langle  v_b \bm{n}_e \cdot \bm{n},\bm{n} \cdot \sigma \cdot \bm{n} \right\rangle_{\partial T} +\left\langle \bm{n} \times \boldsymbol{v}_0,\bm{n} \times\sigma \cdot \bm{n} \right\rangle_{\partial T},
\end{aligned}
\end{equation}
where $\bm n$ represents  the unit outward normal vector to ${\partial T}$ and $ {\bm n}_e$ represents the
 assigned unit normal vector to the  edge or face $e\subset\partial T$. 

\end{definition}

\begin{definition}[modified divergence \cite{su2025novel}]
The modified divergence operator for $\boldsymbol{v} \in \bm V_h $ is defined as  $\nabla_{m} \cdot \boldsymbol{v}\in P_0(\mathcal T_h)$  such that for all $\varphi \in P_0(T)$ and $T\in\mathcal T_h$,
\begin{equation}\label{eq:div}
\left(\nabla_{m} \cdot \boldsymbol{v}, \varphi\right)_T=\left\langle  v_b \bm{n}_e \cdot \bm{n}, \varphi\right\rangle_{\partial T}.
\end{equation}
 
\end{definition}

On the other hand, to develop a pressure-robust EG method,
we introduce a velocity reconstruction operator $\mathcal{R}:  \bm V_h  \rightarrow$ $H(\operatorname{div}; \Omega)$ that maps any function $\bm v \in \bm V_h$ into the lowest-order Raviart-Thomas space  with shape function space $\mathrm{RT}_0=[P_0(T)]^d\oplus \bm x P_0(T)$ such that 
 $$
\int_e \mathcal{R} \bm v \cdot \bm{n}_e \, \d s=\int_e v_b \, \d s \quad \forall e \in \mathcal{E}_h.
$$
Here and hereafter, we denote \[H(\div;\Omega)=\left\{\bm v \in [L^2(\Omega)]^d: \nabla \cdot \bm v \in L^2(\Omega)\right\}.
\]
Denote by $\bm Q_0$ the Scott-Zhang type interpolation operator \cite{scottzhang} from $\left[H^{1/2+\delta}(\Omega)\right]^d$ with $\delta>0$ onto the space of CG and $Q_b$ the $L^2$ projection from $L^2(\mathcal E_h)$ onto $\operatorname{DG}$. 
For $\bm v\in\left[H^{1/2+\delta}(\Omega)\right]^d$, define $v_n\in L^2(\mathcal E_h)$ such that
\[v_n|_e = \bm v|_e\cdot \bm n_e \text{ for any }e\in\mathcal E_h.\]
We define $\bm Q_h\bm v$ by
\[\bm Q_h\bm v =\{\bm Q_0\bm v,Q_b v_n\}.\]
Let $\mathcal{R}_h$ be the lowest-order Raviart-Thomas interpolation, then we have
\begin{align}\label{rh}
\mathcal{R} \bm Q_h {\bm v}=\mathcal{R}_h {\bm v}
.\end{align}
Let $\mathcal Q_h$ and $\mathbb Q_h$ be two $L^2$ projections onto $P_{0}(\mathcal T_h)$ and $\left[P_{0}(\mathcal T_h)\right]^{d \times d}$, respectively.

Define the discrete space with vanishing Dirichlet  boundary conditions
$$
\bm V_h^0=\left\{{\bm v}_h=\left\{{\bm v}_0, v_b\right\} \in \bm V_h, \ {\bm v}_0=0,\ v_b=0 \text { on } \partial \Omega\right\}.
$$
We also define
$$
W_h^0=W_h\cap L_0^2(\Omega),
$$
where $L_0^2(\Omega) = \{v\in L^2(\Omega):(v,1)=0\}.$

We are now in a position to introducing the PR\&PF-EG method for the NS equations with pure vanishing Dirichlet boundary condition, see Algorithm \ref{eq:egn}.

\begin{algorithm}
\caption{A PF\&PR-EG method    for vanishing Dirichlet boundary condition}\label{eq:egn}
\begin{algorithmic}\STATE
    Find ${\bm u}_h =\{\bm u_0,u_b\}\in \bm V_h^0$ and $p_h \in W_h^0$ such that
\begin{align}\label{eqn-9.1.1}
a \left({\bm u}_h,\bm v\right)+ c \left({\bm u}_h, {\bm u}_h,\bm v\right)-b \left(\bm v, p_h\right) & =\left({\bm f}, \mathcal{R} \bm v\right) & & \forall \bm v \in \bm V_h^0, \\ \label{eqn-9.1.2}
b \left({\bm u}_h, q\right) & =0 & & \forall q \in W_h^0,
\end{align}
where 
\begin{align*}
    a(\bm w, \bm v) & :=\nu \left(\nabla_m \bm w, \nabla_m \bm v\right)+s(\bm w, \bm v),\\
s(\bm w, \bm v) & := \nu \sum_{T \in \mathcal{T}_h} h_T^{-1}\left\langle Q_bw_{0,n}-w_b,  Q_bv_{0,n}-v_b\right\rangle_{\partial T}, \\
c(\bm w, \bm z, \bm v) & :=\left( (\nabla\times \bm w_0) \times   \mathcal{R}\bm z,  \mathcal{R}\bm v\right),\\
b(\bm w, q) & :=\left(\nabla_m \cdot \bm w, q\right).
\end{align*}
\end{algorithmic}
\end{algorithm}
\begin{remark}\label{mixboundary}
We also consider the mixed boundary conditions 
$$
\boldsymbol{u}=\boldsymbol{u}_{D} \  \text { on } \Gamma_D\ \text{ and 
 }\ \left(\nu \nabla \boldsymbol{u}-p^{\text {kin }} \mathbb I\right) \boldsymbol{n}=\boldsymbol{u}_N 
 \ \text { on } \Gamma_N,
$$
where  $\Gamma_D$ and $\Gamma_N$  are the Dirichlet  and Neumann boundaries satisfying  $\partial\Omega=\Gamma_D \cup \Gamma_N$ and $\Gamma_D \cap \Gamma_N=\emptyset$,
$\mathbb I$ denotes the $d \times d$ identity matrix, and $\boldsymbol{u}_{D}$, $\boldsymbol{u}_{N}$ are given functions. We define the EG space with boundary conditions
$$
\bm V_h^{0,D}=\left\{{\bm v}_h=\left\{{\bm v}_0, v_b\right\} \in \bm V_h,\ {\bm v}_0=0, \ v_b =0\text { on } \Gamma_D  \right\}.
$$
The PR\&PF-EG method is to find 
${\bm u}_h=\{\bm u_0,u_b\}\in \bm V_h$ with $ \bm u_0=
\Pi_h\bm u_D$, $u_b=Q_b\left(u_{D,n}\right)$ on $\Gamma_D$
and $p_h \in W_h$ 
such that
\begin{align}
a \left({\bm u}_h,\bm v\right)&+ c \left({\bm u}_h, {\bm u}_h,\bm v\right)+d\left({\bm u}_h, {\bm u}_h,\bm v\right)-b \left(\bm v, p_h\right)\nonumber\\
& =\left({\bm f}, \mathcal{R} \bm v\right)
 +\left\langle {\bm n}_e \times {\bm u}_N, {\bm n}_e \times {\bm v}_0\right\rangle_{\Gamma_N}+\left\langle  {\bm u}_N \cdot {\bm n}_e, v_b \right\rangle_{\Gamma_N}
&\forall \bm v \in \bm V_h^{0,D}, \\ 
b \left({\bm u}_h, q\right) & =0 &\forall q \in W_h,
\end{align}
where 
$$
d(\bm u_h, \bm  u_h, \bm v)  := \frac{1}{2} \langle|\bm u_0|^2,v_b \rangle_{\Gamma_N},$$
and $\Pi_h$ represents the first-order Lagrange interpolation. 





\end{remark}


\section{Theoretical analysis} \label{sec:Theor}
In this section, we analyze  the well-posedness and establish error estimates 
for the proposed PR\&PF-EG method in Algorithm \ref{eq:egn}. To this end, we define the following   mesh-dependent norm corresponding to the bilinear form $a(\cdot,\cdot)$ in $\bm V_h^0$:
$$
|\!|\!|\bm v |\!|\!|^2=\sum_{T \in \mathcal{T}_h}\left\|\nabla_m \bm v\right\|_T^2+\sum_{T \in \mathcal{T}_h} h_T^{-1}\left\|Q_bv_{0,n}- v_b\right\|_{\partial T}^2.
$$
\begin{lemma}\label{lem:8}
For any $\boldsymbol{v}=\left\{\boldsymbol{v}_0, v_b\right\} \in \bm V_h^0$, the following inequality holds
\begin{equation}\label{eq:v0}
\left\|\nabla \boldsymbol{v}_0\right\|^2 \leq C |\!|\!| \boldsymbol{v} |\!|\!|^2.
\end{equation}
Therefore, $|\!|\!| \cdot |\!|\!|$ defines a norm in $\bm V_h^0$.
\end{lemma}
\begin{proof}
For any $\bm v=\left\{\bm v_0, v_b\right\} \in \bm V_h^0$, it follows from  
 the integration by parts, \eqref{eq:grad}, and the definition of   $Q_b$ that 
$$
\begin{aligned}
\left(\nabla \bm v_0, \nabla \bm v_0\right)_T  = &-\left(\bm v_0, \nabla \cdot \nabla \bm v_0\right)_T+\left\langle\bm v_0, \nabla \bm v_0 \cdot \bm{n}\right\rangle_{\partial T} \\
 =&\left\langle v_b\bm{n}_e \cdot \bm{n}, \bm{n} \cdot \nabla \bm v_0 \cdot \bm{n}\right\rangle_{\partial T}+ \left\langle \bm{n} \times \bm v_0  ,\bm{n} \times \nabla \bm v_0 \cdot \bm{n}\right\rangle_{\partial T}\\ 
&+\left\langle\bm v_0 \cdot \bm{n} - v_b\bm{n}_e \cdot \bm{n},\bm{n} \cdot \nabla \bm v_0 \cdot \bm{n}\right\rangle_{\partial T} \\
 =&\left(\nabla_m \bm v, \nabla \bm v_0\right)_T+\left\langle (Q_bv_{0,n}-v_b)\bm{n}_e \cdot \bm{n}, \bm{n} \cdot \nabla \bm v_0 \cdot \bm{n}\right\rangle_{\partial T} .
\end{aligned}
$$
By applying the trace inequality and the inverse inequality, we obtain
$$
\left\|\nabla \bm v_0\right\|_T^2 \leq C\left(\left\|\nabla_m \bm v\right\|_T\left\|\nabla \bm v_0\right\|_T+h_T^{-\frac{1}{2}}\left\|Q_bv_{0,n} - v_b\right\|_{\partial T}\left\|\nabla \bm v_0\right\|_T\right),
$$
which implies
$$
\left\|\nabla \bm v_0\right\|_T^2 \leq C\left(\left\|\nabla_m \bm v\right\|_T^2+h_T^{-1}\left\|Q_b v_{0,n}-v_b\right\|_{\partial T}^2\right),
$$
and hence \eqref{eq:v0} after a summation over $T \in \mathcal{T}_h$.

To show that $|\!|\!| \cdot |\!|\!|$ defines a norm in $\bm V_h^0$,  we only verify the positive length property of $|\!|\!| \cdot |\!|\!|$. Assume that $|\!|\!| \bm v |\!|\!|= 0$ for some $\bm v \in \bm V_h^0$. It follows from \eqref{eq:v0} that $\bm v_0=$ const in $\Omega$, which together with $\bm v_0=0$ on $\partial\Omega$ yields $\bm v_0=0$. The assumption $|\!|\!| \bm v |\!|\!|= 0$ also implies $ v_b =Q_bv_{0,n}$, which combined with $\bm v_0=0$ leads to $v_b=0$.
\end{proof}

\subsection{Well-posedness}
It follows from the definition of \( |\!|\!| \cdot |\!|\!| \) and the Cauchy-Schwarz inequality that the following boundedness and coercivity properties hold for the bilinear form \(a(\cdot, \cdot)\).

\begin{lemma} For any $\boldsymbol{v}, \boldsymbol{w} \in \bm V_h^0$, we have
\begin{align}\label{eq:a1}
|a(\boldsymbol{v}, \boldsymbol{w})| & \leq \nu |\!|\!|\boldsymbol{v}|\!|\!| |\!|\!| \boldsymbol{w} |\!|\!|, \\
a(\boldsymbol{v}, \boldsymbol{v}) & =\nu |\!|\!| \boldsymbol{v} |\!|\!|^2 .\label{eq:a2}
\end{align}
\end{lemma}
To show some properties of the trilinear term $c(\cdot,\cdot,\cdot)$, we first present the following lemma.
\begin{lemma}\label{lemma-Rv-bound}It holds for all $r \in[1,6]$ and  $\boldsymbol{v} \in \bm V_h^0$,
\begin{align}\label{eqn-9.2}
\left\|\mathcal{R} \boldsymbol{v} \right\|_{L^r(\Omega)} & \leq C |\!|\!| \boldsymbol{v} |\!|\!|.
\end{align}
Here and hereafter we use $\|\cdot\|_{L^r(\Omega)}$ to denote the $L^r$-norm.
\end{lemma}
\begin{proof}
    By following a similar argument as in the proof of \cite[Proposition 3]{Daniel2020} and utilizing Lemma \ref{le1}, we complete the proof.
\end{proof}
\begin{lemma}  For any $\boldsymbol{v},\boldsymbol{w},\boldsymbol{z} \in \bm V_h^0$, we have
\begin{align}\label{eqn-9.3}
c(\boldsymbol{v}, \boldsymbol{w}, \boldsymbol{w}) & =0, \\
|c(\boldsymbol{v},\boldsymbol{w}, \boldsymbol{z})| & \leq C_{\mathcal N} |\!|\!|\boldsymbol{v}|\!|\!| |\!|\!|\boldsymbol{w}|\!|\!| |\!|\!| \boldsymbol{z} |\!|\!|, \label{eqn-9.4}
\end{align}
where $C_{\mathcal N}$ is a constant independent of $h$.
\end{lemma}

\begin{proof}
  From the definition of $c(\cdot, \cdot, \cdot)$, we have
$$
\begin{aligned}
c(\bm v, \bm w, \bm w) = \left( (\nabla\times\bm v_0) \times \mathcal{R} \bm w, \mathcal{R} \bm w \right)=0.
\end{aligned}
$$
It follows from H\"{o}lder inequality with exponent $(2,4,4)$ and Cauchy-Schwarz inequality that
$$
\begin{aligned}
& \Big |\sum_{T \in \mathcal{T}_h} \int_T (\nabla\times\bm v_0) \times \mathcal{R} \bm w  \cdot \mathcal{R} \bm z\d \bm x\Big| \\
\leq & \sum_{T \in \mathcal{T}_h}\left\|\nabla\times\bm v_0\right\|_{L^2(T)}\left\|\mathcal{R} \bm w\right\|_{L^4(T)}\left\|\mathcal{R} \bm z\right\|_{L^4(T)} \\
\leq &\ C \|\nabla\bm v_0\| \left\|\mathcal{R} \bm w\right\|_{L^4(\Omega)} \left\|\mathcal{R}  \bm z\right\|_{L^4(\Omega)} \\
\leq &\ C_{\mathcal N} |\!|\!|\bm v|\!|\!| |\!|\!|\bm w|\!|\!| |\!|\!|\bm z|\!|\!|,
\end{aligned}
$$
where Lemma \ref{lem:8} and Lemma \ref{lemma-Rv-bound} with $r=4$ have been used to pass to the last line. 
\end{proof}

We present the following (quasi-)commutative properties, which play a crucial role in the subsequent error analysis.

\begin{lemma}[\cite{su2025novel}] 
For $\boldsymbol{v} \in\left[H^1(\Omega)\right]^d$, $\boldsymbol{w} \in H(\operatorname{div}; \Omega)$, $\sigma_h \in\left[P_{0}(\mathcal T_h)\right]^{d \times d}$, and $\varphi_h \in P_{0}(\mathcal T_h)$, the interpolation or projection operators $\bm Q_h, \mathbb {Q}_h$, and $\mathcal Q_h$ satisfy the following (quasi-)commutative properties 
\begin{align}\label{eq:Qhg}
\left(\nabla_m\left(\bm Q_h \boldsymbol{v}\right), \sigma_h \right)_T& =\left(\mathbb Q_h(\nabla \boldsymbol{v}), \sigma_h\right)_T+\left\langle \boldsymbol{n} \times (\bm Q_0 \boldsymbol{v}  - \boldsymbol{v} ),\boldsymbol{n} \times  \sigma_h \cdot \bm{n} \right\rangle_{\partial T}, \\
\left(\nabla_m \cdot\left(\bm Q_h \boldsymbol{w}\right),\varphi_h \right)_T & =\left(\mathcal Q_h(\nabla \cdot \boldsymbol{w}),\varphi_h \right)_T.\label{eq:Qhd}
\end{align}
\end{lemma}
\begin{lemma} [inf-sup condition, \cite{su2025novel}]\label{infsup}
There exists a positive constant $\beta$ independent of $h$ such that
\begin{equation}\label{eq:inf}
\sup _{\boldsymbol{v} \in \bm V_h^0} \frac{b(\boldsymbol{v}, \rho)}{|\!|\!| \boldsymbol{v}|\!|\!| } \geq \beta\|\rho\|,\quad \forall \rho \in W_h^0.
\end{equation}
\end{lemma}
\begin{lemma} \label{le1}
 The operator $\mathcal{R}$ is divergence-preserving, i.e., for all $\boldsymbol{v}=\{\bm v_0,v_b\} \in \bm V_h$, it holds
\begin{equation}\label{div RV}
\nabla \cdot\left(\mathcal{R} \boldsymbol{v}\right)=\nabla_m \cdot \boldsymbol{v}.
\end{equation}
In addition, it holds
$$
\left.\mathcal{R} \boldsymbol{v}\right|_e\cdot\bm n_e = \left.v_b\right|_e\quad\text{ and  }\quad 
\left\|\mathcal{R} \boldsymbol{v}-\boldsymbol{v}_0\right\| \leq C {\color{blue}h}\3bar \bm v_h\3bar.$$
\end{lemma}

\begin{proof}
Using  the integration by parts and the definition of the operator $\mathcal{R} $, we have
\begin{align*}
(\nabla \cdot \mathcal{R} \bm v, q)_{T}&=-(\mathcal{R} \bm v,\nabla q)_{T}+\langle\mathcal{R} \bm v\cdot\bm{n} ,q\rangle_{\partial T}\\
&=\langle v_b\bm n_e\cdot\bm{n} ,q\rangle_{\partial T}=(\nabla_{m} \cdot \bm v,q)_{T}, \ \ \forall q\in P_0(T),
\end{align*}
which implies \eqref{div RV}.

From the fact that $\left.\mathcal{R} \boldsymbol{v}\right|_e\cdot\bm n_e\in P_0(e)$ and the definition of the operator $\mathcal{R} $, we obtain $\left.\mathcal{R} \boldsymbol{v}\right|_e\cdot\bm n_e = \left.v_b\right|_e$.

Any $\bm \phi \in [P_1(\mathcal T_h)]^d$  is uniquely determined by $\int_e \bm \phi\cdot\bm n_eq\d s$ for $q\in P_1(e)$ and $e\in \mathcal E_h$ \cite{monk2003}. 
 Using the standard scaling argument, we have 
\begin{align}\label{RT1}
\left\|\boldsymbol{\phi}\right\|^2 \leq C \sum_{e\in\mathcal E_h} h_e\left\|\boldsymbol{\phi} \cdot \bm{n}_e\right\|^2_e. 
\end{align}
It follows from $\bm v_0 \in \operatorname{CG}$ and $\mathcal{R} \bm v \in \operatorname{RT}_0$ that 
$\left(\mathcal{R} \bm v -\bm v_0\right) \in [P_1(\mathcal T_h)]^d$.
Letting $\bm \phi=\mathcal{R} \bm v -\bm v_0$ in \eqref{RT1}, we derive
\begin{align*}
\|\mathcal{R} \bm v-\bm v_0\|^2 & \leq C \sum_{e\in\mathcal E_h} h_e\left\|\left(\mathcal{R} \bm v-\bm v_0\right) \cdot \bm{n}_e\right\|_e^2 = C \sum_{e\in\mathcal E_h} h_e\left\|\bm v_0 \cdot \bm{n}_e - v_b\right\|_e^2 \\
&\leq C \sum_{e\in\mathcal E_h} h_e\left\|\bm v_0\cdot\bm n_e-Q_bv_{0,n}\right\|_e^2+ C \sum_{e\in\mathcal E_h} h_e\left\|Q_bv_{0,n} - v_b\right\|_e^2,\\
&\leq C \sum_{e\in\mathcal E_h} h_e^3\left|\bm v_0\cdot\bm n_e\right|_{1,e}^2+ C \sum_{e\in\mathcal E_h} h_e\left\|Q_bv_{0,n} - v_b\right\|_e^2,
\end{align*}
which, together with the trace inequality, the inverse inequality, and Lemma \ref{lem:8}, completes the proof.
\end{proof}
\begin{lemma}[\cite{scottzhang,monk2003}]
Assume that $\boldsymbol{w} \in\left[H^{2}(\Omega)\right]^d$ and $\rho \in H^1(\Omega)$. Then we have
\begin{align}\label{eqn-7.1}
&\sum_{T \in \mathcal{T}_h}\|\boldsymbol{w}-\bm Q_0 \boldsymbol{w}\|_{m,T}^2 \leq Ch^{2(l-m)}\|\boldsymbol{w}\|_{l}^2\quad  \text{for }  0 \leq m \leq l \leq 2,\\
& \sum_{T \in \mathcal{T}_h} \left\|\nabla \boldsymbol{w}-\mathbb Q_h(\nabla \boldsymbol{w})\right\|^2_T \leq C h^{2 }\|\boldsymbol{w}\|_{2}^2, \label{eqn-7.2}\\
& \sum_{T \in \mathcal{T}_h}\left\|\rho-\mathcal Q_h \rho\right\|_T^2 \leq C h^{2}\|\rho\|_1^2 ,\label{eqn-7.3}\\
&\sum_{T \in \mathcal{T}_h}\left\| \boldsymbol{w}-\mathcal{R}_h\boldsymbol{w}\right\|^p_{L^p(T)} \leq C h^{p+d-\frac{dp}{2}}\|\boldsymbol{w}\|_1^p\quad \text{for }1\leq p\leq 6.\label{eqn-7.4}
\end{align}
\begin{proof}
    The proof of \eqref{eqn-7.4} follows from a similar argument to that of \cite[Theorem 5.25]{monk2003}.
\end{proof}
\end{lemma}

Having completed all necessary preparations, we are now ready to employ the Leray-Schauder fixed point theorem \cite{layton2008introduction} to establish both the existence and uniqueness of a solution for Algorithm \ref{eq:egn}. To this end, we define a discrete divergence-free subspace \(\bm D_h\) of \(\bm V_h^0\) as:
\[
\bm D_h = \{\bm{v} \in \bm V_h^0 : \nabla_m \cdot \bm{v} = 0\}.
\]
Consequently, the PF\&PR-EG method in Algorithm \ref{eq:egn} is reformulated to find \(\bm{u}_h \in \bm D_h\) satisfying:
\begin{align}\label{eqn-9.8}
a(\bm{u}_h, \bm{v}) + c(\bm{u}_h, \bm{u}_h, \bm{v}) = (\bm{f}, \mathcal{R} \bm{v}), \quad \forall \bm{v} \in \bm D_h.
\end{align}


\begin{lemma}[Stability of the discrete velocity and pressure]
    Any solution $(\bm u_h;p_h)$ of Algorithm \ref{eq:egn} satisfies
    \begin{align}
    |\!|\!| \bm u_h |\!|\!|&\leq \frac{\|{\bm f}\|}{\nu },\label{stability}\\
      \|p_h\|&\leq \left(\frac{2}{\beta }+\frac{C_{\mathcal N}\|{\bm f}\|}{\beta\nu^2}\right)\|{\bm f}\|. \label{stability-p} 
    \end{align} 
\end{lemma}
\begin{proof}
For the proof of \eqref{stability}, set $\bm v=\bm u_h$ in \eqref{eqn-9.8} and use \eqref{eq:a2} -- \eqref{eqn-9.3}. 

    Solve the following for $p_h$:
    \[b \left(\bm v, p_h\right) =-\left({\bm f}, \mathcal{R} \bm v\right)+a \left({\bm u}_h,\bm v\right)+ c \left({\bm u}_h, {\bm u}_h,\bm v\right), \ \forall \bm v\in \bm D_h^{\perp}, \]
    where $\bm D_h^{\perp} = \{\bm v\in \bm V_h^0:a(\bm v,\bm w)=0\text{ for any }\bm w\in \bm D_h\}$.
    Then we have
      \[b \left(\bm v, p_h\right) \leq \|{\bm f}\||\!|\!| \bm v|\!|\!|+\nu|\!|\!|{\bm u}_h|\!|\!||\!|\!|\bm v|\!|\!|+ C_{\mathcal N} |\!|\!|{\bm u}_h|\!|\!|^2|\!|\!|\bm v|\!|\!|, \]
      which leads to
      \[\frac{b \left(\bm v, p_h\right)}{|\!|\!| \bm v|\!|\!|} \leq \|{\bm f}\|+\nu|\!|\!|{\bm u}_h|\!|\!|+ C_{\mathcal N} |\!|\!|{\bm u}_h|\!|\!|^2.\]
      Taking the supremum over $\bm v \in \bm V_h^0$ and using Lemma \ref{infsup} and \eqref{stability} give the bound of discrete pressure. 
\end{proof}
 \begin{theorem}[Existence of the discrete velocity and pressure]\label{th4}
 The PR$\&$PF-EG method in Algorithm \ref{eq:egn} has at least one solution ${\bm u}_h \in \bm D_h$ and $p_h\in W_h^0$. 
\end{theorem}
\begin{proof}

Let $F: \bm D_h \rightarrow \bm D_h$ be a nonlinear map such that for each $\bm w \in \bm D_h, \ \widetilde{{\bm u}}_h:=F(\bm w) \in \bm D_h$ is  the solution of the following linear problem
\begin{align}\label{eqn-9.7}
a\left(\widetilde{{\bm u}}_h, \bm v\right)+c\left(\bm w, \widetilde{{\bm u}}_h, \bm v\right)=\left({\bm f}, \mathcal{R}  \bm v\right), \ \forall \bm v \in \bm D_h.
\end{align}
Proceeding as in the proof of \eqref{stability}, we obtain the uniform boundedness of $\widetilde{\bm u}_h$
\begin{align}\label{uniform-bound}
    |\!|\!| \widetilde{\bm u}_h |\!|\!|\leq \frac{\|{\bm f}\|}{\nu}.
\end{align} 
Letting $\bm f=0$ leads to $\widetilde{\bm u}_h=0$, which implies the uniqueness, and hence the existence of the solution of the finite-dimensional problem \eqref{eqn-9.7}. 
Therefore, the map $F$ is well-defined.

From \eqref{eq:a2}, \eqref{eqn-9.4}, and the uniform boundedness \eqref{uniform-bound}, the map $F$ is continuous and, therefore, compact in the finite dimensional space $\bm D_h$. For $0\leq \lambda\leq 1$, consider \[\lambda F(\bm w_{\lambda})= \bm w_{\lambda}.\]
From the uniform boundedness \eqref{uniform-bound}, we have
$$
|\!|\!| \bm w_{\lambda} |\!|\!|\leq \frac{\lambda\|{\bm f}\|}{\nu }\leq \frac{\|{\bm f}\|}{\nu},
$$
which is uniform with respect to $\lambda$. The Leray-Schauder fixed point theorem implies that the nonlinear map $F$ defined by \eqref{eqn-9.7} has  at least one  fixed point ${\bm u}_h\in \bm D_h$ satisfying 
$$
F\left({\bm u}_h\right)={\bm u}_h.
$$ The fixed point ${\bm u}_h$ is a discrete velocity satisfying \eqref{eqn-9.8}. The existence of discrete pressure $p_h$ follows from the existence of discrete velocity $\bm u_h$ and  the inf-sup condition (Lemma \ref{infsup}). 
\end{proof}

\begin{theorem}[Uniqueness under a small data condition]\label{th5}
Suppose the small data condition
 \[\frac{C_{\mathcal N}\| \boldsymbol{f}\|}{\nu^2}<1,\]
where $C_{\mathcal N}$ is defined in \eqref{eqn-9.4}. 
Then Algorithm \ref{eq:egn} has at most one solution.
\end{theorem}
\begin{proof}
Let ${\bm u}_h$, $\overline{{\bm u}}_h \in \bm D_h$ be two solutions of  \eqref{eqn-9.8} and denote $\bm\phi_h={\bm u}_h-\overline{{\bm u}}_h$, then we have
$$
a\left(\bm\phi_h, \bm v\right)+c\left({\bm u}_h, {\bm u}_h, \bm v\right)-c\left(\overline{{\bm u}}_h, \overline{{\bm u}}_h, \bm v\right)=0,\quad \forall \bm v \in \bm D_h.
$$
Noticing that
$$
c\left({\bm u}_h, {\bm u}_h, \bm v\right)-c\left(\overline{{\bm u}}_h, \overline{{\bm u}}_h, \bm v\right)=c\left(\bm\phi_h, {\bm u}_h, \bm v\right)+c\left(\overline{{\bm u}}_h, \bm\phi_h, \bm v\right), 
$$
we obtain
\begin{align}\label{unique-iden1}
    a\left(\bm\phi_h, \bm v\right)+c\left(\overline{{\bm u}}_h, \bm\phi_h, \bm v\right)=-c\left(\bm\phi_h, {{\bm u}}_h, \bm v\right),\quad \forall \bm v \in \bm D_h.
\end{align}
Letting $\bm v=\bm\phi_h$ in \eqref{unique-iden1}, by \eqref{eqn-9.3} and \eqref{eqn-9.4}, we obtain,
$$
\nu |\!|\!| \bm\phi_h |\!|\!|^2 =\left| c\left(\bm\phi_h, {{\bm u}}_h, \bm\phi_h\right) \right| \leq C_{\mathcal N} |\!|\!|  {{\bm u}}_h |\!|\!| |\!|\!| \bm\phi_h |\!|\!|^2.
$$
Note that ${{\bm u}}_h$ is a solution of \eqref{eqn-9.8}, then applying \eqref{stability} gives
$$
\nu |\!|\!| \bm \phi_h |\!|\!|^2 \leq \frac{C_{\mathcal N}\|{\bm f}\|}{\nu}|\!|\!| \bm \phi_h |\!|\!|^2.
$$
Dividing through by $\nu$, we obtain
\[\left(1 - \frac{C_{\mathcal N}\|{\bm f}\|}{\nu^2}\right)|\!|\!| \bm \phi_h |\!|\!|^2\leq 0,\] 
which yields $\bm\phi_h=0$ since the small data condition implies that the multiplier on the left-hand side is positive.  The uniqueness of discrete velocity, and hence the uniqueness of discrete pressure, then follows. 
\end{proof}

\subsection{Error equations}
We establish error estimates for discrete velocity and pressure with respect to the mesh-dependent norm $|\!|\!| \cdot |\!|\!|$. 
Let ${\bm u}_h=\left\{{\bm u}_0, u_b\right\} \in \bm V_h^0$ and $p_h \in W_h^0$ be  the discrete solution of the PR\&PF-EG method \eqref{eqn-9.1.1}--\eqref{eqn-9.1.2}. Denote by ${\bm u}$ and $p$ the exact solution of NS equations \eqref{eq:n1}.
We recall that 
the projections of ${\bm u}$ and $p$ to the finite element spaces  $\bm V_h$ and $W_h$ are given by
$$
\bm Q_h {\bm u}=\left\{\bm Q_0 {\bm u}, Q_bu_n \right\} \quad\hbox{ and } \quad \mathcal Q_h p,
$$ 
respectively. The corresponding errors ${\bm e}_h$ and $\varepsilon_h$   are defined by
\begin{align}\label{eq:eh}
{\bm e}_h=\left\{{\bm e}_0, e_b\right\}=\left\{\bm Q_0 {\bm u}-{\bm u}_0, Q_bu_n-u_b\right\} \quad \text{ and }\quad\varepsilon_h=\mathcal Q_h p-p_h.
\end{align}


Now we derive the main error equations in the following lemmas.
\begin{lemma}\label{lem:9}
Let $(\boldsymbol{u} ; p) \in\left[H_0^1(\Omega)\cap H^2(\Omega)\right]^d \times [H^1(\Omega)\cap L_0^2(\Omega)]$ satisfy  equations \eqref{eq:n1}. Then we have
\begin{equation}\label{eq:EQ}
\nu \left(\nabla_m\left(\bm Q_h \bm u\right), \nabla_m \boldsymbol{v}\right)-\left(\nabla_m \cdot \boldsymbol{v}, \mathcal Q_h p\right)=\left(\boldsymbol{f}, \mathcal{R} \boldsymbol{v}\right)+\ell_{\bm u}(\boldsymbol{v})+\theta_{\bm u}(\boldsymbol{v})+\varphi_{\bm u}(\bm v),\ \forall \boldsymbol{v} \in \bm V_h^0,
\end{equation}
 where $\ell_{\bm u}(\boldsymbol{v})$, $\theta_{\bm u}(\boldsymbol{v})$, and $\varphi_{\bm u}(\bm v)$ are defined by
$$
\begin{aligned}
\ell_{\bm u}(\boldsymbol{v}) =&\nu \sum_{T \in \mathcal{T}_h}\left\langle\boldsymbol{v}_0\cdot\boldsymbol{n}-v_b\boldsymbol{n}_e \cdot \boldsymbol{n}, \boldsymbol{n} \cdot \left(\nabla \bm u \cdot \boldsymbol{n}-\mathbb Q_h(\nabla \bm u) \cdot \boldsymbol{n}\right)\right\rangle_{\partial T} \\
&+\nu \sum_{T \in \mathcal{T}_h}\left\langle \bm{n} \times (\bm Q_0 \bm u  - \bm u) , \bm{n} \times \nabla_m \boldsymbol{v} \cdot \boldsymbol{n}\right\rangle_{\partial T},\\
\theta_{\bm u}(\boldsymbol{v})&=-\nu(\Delta \bm u, \boldsymbol{v}_0-\mathcal{R} \boldsymbol{v}),\\
\varphi_{\bm u}(\bm v)&=-\left((\nabla\times\bm u) \times  \bm u, \mathcal{R} \boldsymbol{v}\right).
\end{aligned}
$$
\end{lemma}
\begin{proof} First, it follows from \eqref{eq:grad}, \eqref{eq:Qhg}, and the integration by parts that
\begin{equation}\label{eqn-11}
\begin{aligned}
\left(\nabla_m\left(\bm Q_h \bm u\right), \nabla_m \bm v\right)_T  =&\ \left(\mathbb Q_h(\nabla \bm u), \nabla_m \bm v\right)_T +\left\langle \boldsymbol{n} \times (\bm Q_0 \bm u  - \bm u ),\boldsymbol{n} \times \nabla_m \bm v \cdot \bm{n}\right\rangle_{\partial T}\\
 =&\ -\left(\bm v_0, \nabla \cdot \mathbb Q_h(\nabla \bm u)\right)_T+\left\langle v_b\bm{n}_e \cdot \bm{n},\bm{n} \cdot \mathbb Q_h(\nabla \bm u) \cdot \bm{n}\right\rangle_{\partial T} \\
& +\left\langle\boldsymbol{n} \times \bm v_0,\boldsymbol{n} \times \mathbb Q_h(\nabla \bm u) \cdot \bm{n}\right\rangle_{\partial T}\\
&+\left\langle \boldsymbol{n} \times (\bm Q_0 \bm u  - \bm u) , \boldsymbol{n} \times \nabla_m \bm v \cdot \bm{n}\right\rangle_{\partial T}\\
 =&\ \left(\nabla \bm v_0, \mathbb Q_h(\nabla \bm u)\right)_T-\left\langle \bm v_0\cdot\bm{n}-v_b\bm{n}_e \cdot \bm{n},\bm{n} \cdot \mathbb Q_h(\nabla \bm u) \cdot \bm{n}\right\rangle_{\partial T} \\
& +\left\langle \boldsymbol{n} \times (\bm Q_0 \bm u  - \bm u) , \boldsymbol{n} \times \nabla_m \bm v \cdot \bm{n}\right\rangle_{\partial T}\\
 =&\ \left(\nabla \bm u, \nabla \bm v_0\right)_T-\left\langle\bm v_0\cdot\bm{n}-v_b\bm{n}_e \cdot \bm{n},\bm{n} \cdot \mathbb Q_h(\nabla \bm u) \cdot \bm{n}\right\rangle_{\partial T}\\
& +\left\langle \boldsymbol{n} \times (\bm Q_0 \bm u  - \bm u) , \boldsymbol{n} \times \nabla_m \bm v \cdot \bm{n}\right\rangle_{\partial T}.
\end{aligned}
\end{equation}
By using Lemma \ref{le1}, the integration by parts, and the facts that $\mathcal{R} \bm v \cdot \bm{n}_e$ is continuous across two elements,  we obtain
\begin{equation}\label{eq:b1}
 \left(\nabla_m \cdot \boldsymbol{v}, \mathcal Q_h p\right)=(\nabla \cdot \mathcal{R} \bm v,p)=-(\nabla p, \mathcal{R} \bm v).   
\end{equation}
We multiply \eqref{eq:n1a} by $\mathcal{R} \bm v$ and integrate over $\Omega$
to get
\begin{align}\label{eqn-12}
\left({\bm f},\mathcal{R} \bm v \right)&=-\nu \left(\Delta \bm u, \mathcal{R} \bm v\right)+
\left((\nabla\times\bm u) \times \bm u,\mathcal{R} \bm v\right)+\left(\nabla p, \mathcal{R} \bm v\right)\nonumber\\
&=-\nu \left(\Delta \bm u, \bm v_0\right)+\nu \left(\Delta \bm u, \bm v_0-\mathcal{R} \bm v\right)+
\left((\nabla\times\bm u) \times \bm u,\mathcal{R} \bm v\right)-\left(\nabla_m \cdot \boldsymbol{v}, \mathcal Q_h p\right),
\end{align}
where we have used \eqref{eq:b1}. 
It follows from the integration by parts that
$$
-\nu \left(\Delta \bm u, \bm v_0\right)=\nu \sum_{T \in \mathcal{T}_h}\left(\nabla \bm u, \nabla \bm v_0\right)_T,
$$
which, combined with \eqref{eqn-11}, gives
\begin{equation}\label{eqn-13}
\begin{aligned}
-\nu \left(\Delta \bm u, \bm v_0\right)& =\nu \left(\nabla_m\left(\bm Q_h \bm u\right), \nabla_m \bm v\right) \\
& -\nu\sum_{T \in \mathcal{T}_h}\left\langle \bm v_0 \cdot \bm{n}-v_b\bm{n}_e \cdot \bm{n},\bm{n} \cdot \left(\nabla \bm u \cdot \bm{n}-\mathbb Q_h(\nabla \bm u) \cdot \bm{n}\right)\right\rangle_{\partial T}\\
& -\nu\sum_{T \in \mathcal{T}_h}\left\langle  \bm{n} \times (\bm Q_0 \bm u  - \bm u) ,\bm{n} \times \nabla_m \bm v \cdot \bm{n}\right\rangle_{\partial T},
\end{aligned}
\end{equation}
where we have used the fact that $\sum_{T \in \mathcal{T}_h}\left\langle v_b \bm{n}_e \cdot \bm{n}-\bm v_0\cdot\bm n,\bm{n} \cdot  \nabla \bm u \cdot \bm{n}\right\rangle_{\partial T}=0$.
Substituting \eqref{eqn-13} into \eqref{eqn-12} yields \eqref{eq:EQ},
which completes the proof.
\end{proof}

\begin{lemma}\label{err-eqn}
For any $\boldsymbol{v} \in \bm V_h^0$ and $q \in W_h^0$, we have
\begin{align}\label{eqn-9.5}
 a\left(\boldsymbol{e}_h, \boldsymbol{v}\right)+c\left(\boldsymbol{e}_h, \bm Q_h \boldsymbol{u},\boldsymbol{v} \right)+c\left(\boldsymbol{u}_h, \boldsymbol{e}_h, \boldsymbol{v} \right)-b\left(\boldsymbol{v}, \varepsilon_h\right)& =\chi_{\boldsymbol{u}}(\boldsymbol{v})+\phi_{\boldsymbol{u}}(\boldsymbol{v}), \\
 b\left({\bm e}_h, q\right) & =0,\label{eqn-9.6}
\end{align}
where
$$
\begin{aligned}
\phi_{\boldsymbol{u}}(\boldsymbol{v})  & = c\left(\bm Q_h \boldsymbol{u}, \bm Q_h \boldsymbol{u},\boldsymbol{v} \right)-\left((\nabla\times\boldsymbol{u}) \times  \boldsymbol{u}, \mathcal{R} \boldsymbol{v}\right),\\
\chi_{\boldsymbol{u}}(\boldsymbol{v}) &  = \ell_{\boldsymbol{u}} (\boldsymbol{v})+s\left(\bm Q_h \boldsymbol{u}, \boldsymbol{v}\right)+\theta_{\boldsymbol{u}}(\boldsymbol{v}).
\end{aligned}
$$
\end{lemma}

\begin{proof}
It follows from Lemma \ref{lem:9} that
$$
\nu \left(\nabla_m\left(Q_h {\bm u}\right), \nabla_m \bm v\right)-\left(\nabla_m \cdot \bm v, \mathcal Q_h p \right)=\left({\bm f}, \mathcal{R}\bm v\right)+\ell_{{\bm u}}(\bm v)+\theta_{{\bm u}}(\bm v)-\left((\nabla\times\bm u) \times {\bm u},\mathcal{R} \bm v\right).
$$
Adding $s\left(\bm Q_h {\bm u}, \bm v\right)$ and $c\left(\bm Q_h {\bm u}, \bm Q_h {\bm u}, \bm v\right)$ to both sides of the above equations, we obtain
$$
a \left(\bm Q_h {\bm u}, \bm v\right)+c\left(\bm Q_h {\bm u}, \bm Q_h {\bm u}, \bm v\right)-b \left(\bm v, \mathcal Q_h p\right)=({\bm f}, \mathcal{R} \bm v)+\chi_{{\bm u}}(\bm v)+\phi_{{\bm u}}(\bm v)
.$$
Subtracting \eqref{eqn-9.1.1} from above equation and noticing that
$$
\begin{aligned}
c\left(\bm Q_h {\bm u},\bm Q_h {\bm u}, \bm v\right)-c\left({\bm u}_h, {\bm u}_h, \bm v\right)
=c\left({\bm e}_h, \bm Q_h {\bm u}, \bm v\right)+c\left({\bm u}_h, {\bm e}_h, \bm v\right),
\end{aligned}
$$
yields \eqref{eqn-9.5}.

Multiplying \eqref{eq:n1b} by $q \in W_h^0$ and using \eqref{eq:Qhd} gives
\begin{align*}
0=(\nabla \cdot {\bm u}, q)=\left(\nabla_m \cdot \bm Q_h {\bm u}, q\right).
\end{align*}
Subtracting \eqref{eqn-9.1.2} from above equation yields \eqref{eqn-9.6}.
\end{proof}
\begin{lemma}\label{eq:chi}
Assume that $\bm{w} \in\left[H^2(\Omega)\right]^d$, we have   
\begin{align}
\left|\chi_{\bm{w}}(\bm{v})\right| \leq C \nu h \|\bm{w}\|_{2} |\!|\!|\bm{v} |\!|\!|,\quad \forall \bm{v} \in \bm V_h.
\end{align}
\end{lemma}
\begin{proof}
 For $s\left(\bm Q_h \bm w, \bm v\right)$,  we use the definition of $Q_b$, the trace inequality, and \eqref{eqn-7.1} to have
\begin{align}
\left|s\left(\bm Q_h \bm w, \bm v\right)\right| & =\nu \Big|\sum_{T \in \mathcal{T}_h} h_T^{-1}\left\langle Q_b(Q_0w)_n-Q_b w_n , Q_b v_{0,n}-v_b\right\rangle_{\partial T}\Big| \nonumber\\
& =\nu \Big|\sum_{T \in \mathcal{T}_h} h_T^{-1}\left\langle (Q_0w)_n- w_n , Q_b v_{0,n}-v_b\right\rangle_{\partial T}\Big| \nonumber\\
& \leq \nu \Big(\sum_{T \in \mathcal{T}_h}\left(h_T^{-2}\left\|\bm Q_0 \bm w - \bm w\right\|_T^2+\left\|\nabla\left(\bm Q_0 \bm w-\bm w\right)\right\|_T^2\right)\Big)^{1/2}|\!|\!|\bm v|\!|\!|\nonumber\\
& \leq C \nu h  \|\bm w\|_{2}|\!|\!|\bm v|\!|\!|.\label{eq:s}
\end{align}
  For $\ell_{\bm w}(\bm v)$, we first use  the trace inequality, the inverse inequality, \eqref{eqn-7.2}, and Lemma \ref{le1} to obtain
$$
\begin{aligned}
&\Big|\sum_{T \in \mathcal{T}_h}\left\langle\bm v_0 \cdot \bm{n}-v_b \bm{n}_e \cdot \bm{n}, \bm{n}\cdot (\nabla \bm w \cdot \bm{n}-\mathbb Q_h(\nabla \bm w) \cdot \bm{n})\right\rangle_{\partial T}\Big| \\
\leq &\ C \Big(\sum_{T \in \mathcal{T}_h} h_T^{-1}\left\| \nabla \bm w -\mathbb Q_h(\nabla \bm w) \right\|_{\partial T}^2\Big)^{1 / 2}\Big(\sum_{T \in \mathcal{T}_h} h_T\|\bm v_0 \cdot \bm{n}-v_b \bm{n}_e \cdot \bm{n}\|_{\partial T}^2\Big)^{1/2}\\
\leq &\ C \|\bm w\|_{2}\Big(\sum_{T \in \mathcal{T}_h} h_T\|\bm v_0 -\mathcal R\bm v\|_{\partial T}^2\Big)^{1/2}
\leq \ C \|\bm w\|_{2}\|\bm v_0-\mathcal R\bm v\|\leq C h\|\bm w\|_{2}|\!|\!|\bm v|\!|\!|.\end{aligned}
$$
It follows from the trace inequality  and \eqref{eqn-7.1} that
$$
\begin{aligned}
\Big|\sum_{T \in \mathcal{T}_h}\left\langle \bm{n} \times (\bm Q_0 \bm w  - \bm w) , \bm{n} \times \nabla_m \bm v \cdot \bm{n}\right\rangle_{\partial T}\Big|
\leq  C h \|\bm w\|_{2}|\!|\!|\bm v|\!|\!|.
\end{aligned}
$$
Therefore, we have
\begin{equation*}\label{eqn-6.2}
{\left|\ell_{\bm w}(\bm v)\right| \leq  C \nu h  \|\bm w\|_{2}|\!|\!|\bm v|\!|\!|.}
\end{equation*}
For $\theta_{\bm w}(\bm v)$, we use Lemma \ref{le1}  to obtain
$$
\left|\theta_{\bm w}(\bm v)\right| \leq \nu\|\Delta \bm w\|\|\bm v_0-\mathcal{R} \bm v\|\leq  C \nu h\|\bm w\|_2|\!|\!| \bm v|\!|\!|.
$$
This completes the proof.
\end{proof}
\begin{lemma}\label{eq:phi} Let $\boldsymbol{w} \in\left[H_0^1(\Omega)\cap H^2 (\Omega)\right]^d$ and $\boldsymbol{v} \in \bm V_h^0$, then we have
\begin{align}
\left|\phi_{\boldsymbol{w}}(\boldsymbol{v})\right| \leq C h \|\nabla \boldsymbol{w}\| \|\boldsymbol{w}\|_2 |\!|\!|\boldsymbol{v} |\!|\!| .
\end{align}
\end{lemma}

\begin{proof}
For any $\bm v \in \bm V_h^0$, by the definition of $c(\cdot,\cdot,\cdot)$, we have
$$
\begin{aligned}
& c\left(\bm Q_h  \bm w,\bm Q_h  \bm w,\bm v \right)-\left((\nabla\times \bm w) \times  \bm w, \mathcal{R}  \bm v\right)\\
= & \sum_{T \in \mathcal{T}_h} \int_T  \left(\nabla\times \bm Q_0 \bm w \right) \times \mathcal{R} \left(\bm Q_h \bm w \right)\cdot \mathcal{R} \bm v  -  (\nabla\times \bm w)  \times  \bm w \cdot \mathcal{R} \bm v \d \bm x\\
= & \sum_{T \in \mathcal{T}_h} \int_T  \nabla\times \left( \bm Q_0 \bm w - \bm w\right)  \times \left(\mathcal{R} \left(\bm Q_h \bm w \right)-\bm w\right)\cdot \mathcal{R} \bm v \d \bm x\\
&+ \sum_{T \in \mathcal{T}_h} \int_T  \nabla\times \left( \bm Q_0 \bm w - \bm w \right)  \times \bm w \cdot \mathcal{R} \bm v \d \bm x\\
&+ \sum_{T \in \mathcal{T}_h} \int_T  (\nabla\times  \bm w )  \times \left(\mathcal{R} \left(\bm Q_h \bm w \right)-\bm w\right)\cdot \mathcal{R} \bm v \d \bm x\\
=:& \operatorname{I}+\operatorname{II}+\operatorname{III}.
\end{aligned}
$$
By H\"{o}lder inequalities with exponents (2, 4, 4), Cauchy-Schwarz inequality, \eqref{rh}, \eqref{eqn-7.1}, \eqref{eqn-7.4}, Poincar\'{e} inequality, and \eqref{eqn-9.2}, we have
$$
\begin{aligned}
\operatorname{I}&=  \sum_{T \in \mathcal{T}_h} \int_T  \nabla\times \left( \bm Q_0 \bm w - \bm w\right)  \times \left(\mathcal{R} \left(\bm Q_h \bm w \right)-\bm w\right)\cdot \mathcal{R} \bm v \d \bm x \\
&\leq  \sum_{T \in \mathcal{T}_h}\left\|\nabla\times \left( \bm Q_0 \bm w - \bm w\right)\right\|_{L^2(T)}\left\|\mathcal{R} \left(\bm Q_h \bm w \right)-\bm w \right\|_{L^4(T)}\left\|\mathcal{R} \bm v\right\|_{L^4(T)} \\
&\leq  C\left\|\nabla \left( \bm Q_0 \bm w - \bm w\right)\right\|_{L^2(\Omega)} \left\|\mathcal{R}_h \bm w -\bm w\right\|_{L^4(\Omega)} \left\|\mathcal{R} \bm v\right\|_{L^4(\Omega)} \\
&\leq Ch\|\bm w\|_2 \|\nabla \bm w\|  |\!|\!| \bm v |\!|\!|.
\end{aligned}
$$
By H\"{o}lder inequalities with exponents (2, 4, 4), Cauchy-Schwarz inequality, the embedding  $H^1(\Omega) \hookrightarrow L^4(\Omega)$, \eqref{eqn-7.1}, Poincar\'{e} inequality, and \eqref{eqn-9.2}, we have
$$
\begin{aligned}
\operatorname{II}&= \sum_{T \in \mathcal{T}_h} \int_T  \nabla\times \left( \bm Q_0 \bm w - \bm w \right)  \times \bm w \cdot \mathcal{R} \bm v \d \bm x\\
&\leq  C\left\|\nabla \left( \bm Q_0 \bm w - \bm w\right)\right\|_{L^2(\Omega)}\left\|\bm w\right\|_{L^4(\Omega)}\left\|\mathcal{R} \bm v\right\|_{L^4(\Omega)} \\
&\leq  Ch\|\bm w\|_2\|\nabla\bm w\||\!|\!| \bm v |\!|\!|.
\end{aligned}
$$
It follows from \eqref{rh}, \eqref{eqn-7.4}, the embedding  $H^1(\Omega) \hookrightarrow L^4(\Omega)$, Poincar\'{e} inequality, and \eqref{eqn-9.2} that
$$
\begin{aligned}
\operatorname{III}&= \sum_{T \in \mathcal{T}_h} \int_T  (\nabla\times  \bm w)  \times  \left(\mathcal{R} \left(\bm Q_h \bm w \right)-\bm w\right) \cdot \mathcal{R} \bm v \d \bm x \\
&\leq  C\left\|\nabla \bm w\right\|_{L^4(\Omega)}\left\|\mathcal{R}_h \bm w -\bm w\right\|_{L^2(\Omega)}\left\|\mathcal{R} \bm v\right\|_{L^4(\Omega)} \\
&\leq  Ch\|\nabla\bm w\|_1\|\nabla\bm w\|
|\!|\!| \bm v |\!|\!|.
\end{aligned}
$$
Combining the estimates of I, II, and III above completes the proof.
\end{proof}

\subsection{Error estimates}
According to Helmholtz decomposition, we decompose ${\bm f}=\bm g+\nabla \psi$, where $\bm g\in H(\div;\Omega)$ such that $\nabla\cdot\bm g=0$ and $\psi \in H^1(\Omega)$. We multiply \eqref{eq:n1a} by $\bm u$, integrate over $\Omega$, and use \eqref{eq:n1b}, \eqref{eq:n1c}, and $((\nabla\times\bm u)\times\bm u,\bm u)=0$ to obtain
$$
\nu\|\nabla {\bm u}\|^2=(\bm f, {\bm u}) =(\bm g, {\bm u}) +(\nabla\psi,\bm u) =  (\bm g, {\bm u})\leq\|\bm g\|\|{\bm u}\| \leq C_P\|\bm g\|\|\nabla {\bm u}\|,
$$
where we have used Poincar\'e inequality in the last inequality. Thus, the exact solution is bounded by the solenoidal part of ${\bm f}$ :
\begin{align}\label{eq:ug}
 \|\nabla {\bm u}\| \leq C_P \nu^{-1}\|\bm g\|.  
\end{align}
Similarly, taking $\bm v=\bm u_h$ in \eqref{eqn-9.8} and using Lemmas \ref{lemma-Rv-bound} and \ref{le1} yield
\begin{align}\label{eq:uhg}
 \3bar {\bm u}_h\3bar \leq C \nu^{-1}\|\bm g\|.  
\end{align}
By setting $\sigma_h = \nabla_m \bm Q_h {\bm u}$ in \eqref{eq:Qhg}, we derive
\begin{align*}
    \sum_{T \in \mathcal{T}_h}\left\|\nabla_m \bm Q_h {\bm u}\right\|_T^2 \leq  \|\nabla\bm u\|^2 + C\sum_{T\in\mathcal T_h}h_T^{-1}\|\bm Q_0\bm u-\bm u\|_{\partial T}^2,
\end{align*}
and therefore, 
\begin{align*}
    |\!|\!|\bm Q_h {\bm u} |\!|\!|^2  &=  \sum_{T \in \mathcal{T}_h}\left\|\nabla_m \bm Q_h {\bm u}\right\|_T^2+\sum_{T \in \mathcal{T}_h} h_T^{-1}\left\|Q_b (Q_0u)_n- Q_b(u_n)\right\|_{\partial T}^2 \nonumber\\
    &\leq \sum_{T \in \mathcal{T}_h}\left\|\nabla_m \bm Q_h {\bm u}\right\|_T^2+\sum_{T \in \mathcal{T}_h} h_T^{-1}\left\|\bm Q_0\bm u- \bm u\right\|_{\partial T}^2\nonumber\\
    &\leq \|\nabla\bm u\|^2 + C\sum_{T\in\mathcal T_h}h_T^{-1}\|\bm Q_0\bm u-\bm u\|_{\partial T}^2\nonumber\\
    &\leq C_b^2\|\nabla\bm u\|^2,
\end{align*}
which, together with \eqref{eq:ug}, leads to
\begin{align}\label{bound-Qhu}
    |\!|\!|\bm Q_h {\bm u} |\!|\!|\leq C_bC_P\nu^{-1}\|\bm g\|.
\end{align}

Now we present the error estimates of the proposed PR\&PF-EG method in Algorithm \ref{eq:egn}. 
\begin{theorem}\label{th:error}
Assume that 
\begin{align}\label{eq:g}
  \|\boldsymbol{g}\| \leq \frac{\nu^2}{2C_{\mathcal N}C_bC_P}.
\end{align}
Let $(\boldsymbol{u}; p) \in\left[H_0^1(\Omega)\cap H^2(\Omega)\right]^d\times [H^1(\Omega)\cap L_0^2(\Omega)]$  and $\left(\boldsymbol{u}_h; p_h\right) \in \bm V_h^0 \times W_h^0$ be the solutions of \eqref{eq:n1} and Algorithm \ref{eq:egn}, respectively. Then the following error estimates hold
\begin{align}\label{estimateforv}
\|\bm e_0\|_1&\leq C|\!|\!|\boldsymbol{e}_h |\!|\!|  \leq C h \|\boldsymbol{u}\|_2, \\\label{estimateforp}
\left\|\varepsilon_h\right\| & \leq C\nu h \|{\bm u}\|_2.
\end{align}
\end{theorem}
\begin{proof}
To estimate ${\bm e}_h$, we let $q=\varepsilon_h$ and $\bm v={\bm e}_h$ in the error equations \eqref{eqn-9.5} and \eqref{eqn-9.6} to get
$$
a\left({\bm e}_h, {\bm e}_h\right)+c\left({\bm e}_h,\bm Q_h {\bm u}, \bm e_h\right)=\chi_{{\bm u}}({\bm e}_h)+\phi_{{\bm u}}\left({\bm e}_h\right) .
$$
It follows from   \eqref{eqn-9.4} and \eqref{bound-Qhu} that
$$
\begin{aligned}
\text {LHS} & :=a\left({\bm e}_h, {\bm e}_h\right)+c\left({\bm e}_h, \bm Q_h {\bm u}, {\bm e}_h\right)=\nu |\!|\!| {\bm e}_h |\!|\!|^2+c\left({\bm e}_h, \bm Q_h {\bm u}, {\bm e}_h\right) \\
& \geq\left(\nu-C_{\mathcal N}|\!|\!| \bm Q_h {\bm u} |\!|\!| \right) |\!|\!| {\bm e}_h |\!|\!| ^2 \\
& \geq\left(\nu-\nu^{-1} C_{\mathcal N} C_bC_P\|\bm g\|\right)|\!|\!| {\bm e}_h |\!|\!|^2 \\
& \geq \frac{\nu}{2} |\!|\!| {\bm e}_h |\!|\!| ^2,
\end{aligned}
$$
where we have used \eqref{eq:g}. 
By using Lemmas \ref{eq:chi} -- \ref{eq:phi}, \eqref{eq:ug}, and \eqref{eq:g}, we obtain
$$
\text {RHS}:=\chi_{{\bm u}}({\bm e}_h)+\phi_{{\bm u}}\left({\bm e}_h\right) \leq C \nu h \|{\bm u}\|_2|\!|\!|{\bm e}_h|\!|\!|+C h \|\nabla{\bm u}\| \|{\bm u}\|_{2}|\!|\!|{\bm e}_h|\!|\!|\leq C \nu h \|{\bm u}\|_2|\!|\!|{\bm e}_h|\!|\!|.
$$
Combining the estimates of LHS and RHS and using Lemma \ref{lem:8} above gives \eqref{estimateforv}.

To estimate $\varepsilon_h$, we use Lemmas \ref{infsup} and \ref{err-eqn}, \eqref{eq:a1}, \eqref{eq:uhg}, and \eqref{bound-Qhu} to have
$$
\begin{aligned}
\left\|\varepsilon_h\right\| & \leq \sup _{\bm v \in \bm V_h^0} \frac{\left|b\left(\bm v, \varepsilon_h\right)\right|}{|\!|\!|\bm v|\!|\!|}\\
& =\sup _{\bm v \in \bm V_h^0} \frac{\big| a\left({\bm e}_h, \bm v\right)+c\left({\bm e}_h, \bm Q_h {\bm u}, \bm v\right)+c\left({\bm u}_h,  {\bm e}_h, \bm v\right)-\chi_{{\bm u}}(\bm v)-\phi_{{\bm u}}(\bm v)\big|}{|\!|\!| \bm v|\!|\!|} \\
& \leq  \nu |\!|\!| {\bm e}_h |\!|\!|+ C_{\mathcal N} |\!|\!| {\bm e}_h |\!|\!| |\!|\!|\bm Q_h {\bm u}|\!|\!|+ C_{\mathcal N} |\!|\!| {\bm e}_h |\!|\!| |\!|\!|{\bm u}_h|\!|\!|+ C \nu h \|{\bm u}\|_{2} \\
& \leq C\nu h  \|{\bm u}\|_2.
\end{aligned}
$$
\end{proof}

 \section{Numerical experiments}\label{sec5} 

In this section, we present several numerical examples  to validate  the convergence and demonstrate the   pressure-robustness of the proposed method.
We employ the standard Newton’s iteration method  to  linearize \eqref{eqn-9.1.1}--\eqref{eqn-9.1.2}. To be specific, 
the iteration algorithm is given by: 
For  $n=0,1,...$, we solve  the linearized discrete problem:
$$
\!a\!\left(\bm u_h^{n+1}, \bm v^{n+1}\right)+c\!\left(\bm u_h^n, \bm u_h^{n+1}, \bm v\right)+c\!\left(\bm u_h^{n+1}, \bm u_h^n, \bm v\right)+b\!\left(\bm v, p_h^{n+1}\right)\!=\!\left(\bm{f}, \mathcal {R}\bm v\right)+c\!\left(\bm u_h^n, \bm u_h^n, \bm v\right)\!,
$$
where we take $\left(\bm u_h^0; p_h^0\right)=\bm 0$ unless otherwise specified, 
and use 
$\frac{\left\|\left(\bm u_h^{n+1}; p_h^{n+1}\right)-\left(\bm u_h^n; p_h^n\right)\right\|}{\left\|\left(\bm u_h^{n+1}; p_h^{n+1}\right) \right\|}$ $< 10^{-7} $ or  $n >1000$ 
as the stopping criterion.

\begin{example}[Accuracy test in 2D]\label{Vertex}  \rm
In this example, we confirm the optimal convergence
orders of our proposed  method by considering the following velocity field $\bm u$ and pressure $p$  on $\Omega=(0,1) ^2$:
$$
{\bm u}=\left(\begin{array}{c}
10 x^2(x-1)^2 y(y-1)(2 y-1) \\
-10 x(x-1)(2 x-1) y^2(y-1)^2
\end{array}\right), \quad p=10(2 x-1)(2 y-1).
$$
We impose the pure Dirichlet boundary condition. The body force ${\bm f}$ and boundary velocity field  $\bm u_D = \bm u|_{\partial \Omega}$  are  obtained from the exact solutions.

We implement the proposed method on a uniform triangular mesh with varying viscosities, including $\nu = 1$ and $ \nu =10^{-5}$.  Numerical results are listed in Table \ref{vortex_table},  which verifies the  optimal convergence orders. Comparing the results    of the cases $\nu=1$ and $\nu=10^{-5}$, the error magnitudes remain  almost constant, which validates the pessure-robustness of the method.

We observe that when \(\nu = 10^{-5}\), Newton's iteration fails to converge with an initial guess of zero. To address this, we initially solve the problem for \(\nu = 10^{-3}\) using a zero initial guess. Subsequently, we progressively halve \(\nu\) until reaching \(1.5625 \times 10^{-5}\), after which we directly set \(\nu\) to \(10^{-5}\) for the final step. At each step, the solution from the previous step is used as the initial guess for the next. This gradual reduction in \(\nu\) help in facilitating convergence at the lower viscosity levels.

\begin{table}[h]
\centering
\caption{Example \ref{Vertex}: Numerical results of the proposed method with $\nu=1$ and $\nu=10^{-5}$.}
\label{vortex_table}
\begin{tabular}{@{}lcccccc@{}}
\toprule
$h$ & $\left\|{\bm u}-{\bm u}_0\right\|$ & order & $|{\bm u}-{\bm u}_0|_1$ & order & $\left\|p-p_h\right\|$ & order \\
\midrule
\multicolumn{7}{c}{$\nu=1$} \\
\midrule
1/16 & $1.440\mathrm{e}-03$ & & $8.004\mathrm{e}-02$ & & $3.402\mathrm{e}-01$ & \\
1/32 & $3.640\mathrm{e}-04$ & 1.98 & $4.026\mathrm{e}-02$ & 0.99 & $1.702\mathrm{e}-01$ & 1.00 \\
1/64 & $9.134\mathrm{e}-05$ & 1.99 & $2.017\mathrm{e}-02$ & 1.00 & $8.509\mathrm{e}-02$ & 1.00 \\
1/128 & $2.287\mathrm{e}-05$ & 2.00 & $1.009\mathrm{e}-02$ & 1.00 & $4.254\mathrm{e}-02$ & 1.00 \\
\midrule
\multicolumn{7}{c}{$\nu=10^{-5}$} \\
\midrule
1/16 & $1.659\mathrm{e}-03$ & & $9.512\mathrm{e}-02$ & & $3.400\mathrm{e}-01$ & \\
1/32 & $4.222\mathrm{e}-04$ & 1.97 & $4.189\mathrm{e}-02$ & 1.18 & $1.701\mathrm{e}-01$ & 1.00 \\
1/64 & $1.077\mathrm{e}-04$ & 1.97 & $2.031\mathrm{e}-02$ & 1.04 & $8.504\mathrm{e}-02$ & 1.00 \\
1/128 & $2.720\mathrm{e}-05$ & 1.98 & $1.011\mathrm{e}-02$ & 1.01 & $4.252\mathrm{e}-02$ & 1.00 \\
\bottomrule
\end{tabular}
\end{table}

\end{example}

\begin{example}[Accuracy test in 3D]\label{ex:3D}\rm
In this example, we 
    consider  a 3D problem on the domain  $(0,1)^3$ where   the velocity field $\bm u$ and pressure $p$ are given by \cite{lu2022stabilized}  $$
{\bm u}=\left(\begin{array}{cc}
\alpha \zeta(x)\left(\zeta^{\prime}(y) \zeta(z)-\zeta(y) \zeta^{\prime}(z)\right) \\
\alpha \zeta(y)\left(-\zeta^{\prime}(x) \zeta(z)+\zeta(x) \zeta^{\prime}(z)\right) \\
\alpha \zeta(z)\left(\zeta^{\prime}(x) \zeta(y)-\zeta(x) \zeta^{\prime}(y)\right)
\end{array}\right), \quad p=\alpha \cos (\pi x) \cos (\pi y) \sin (\pi z).
$$
Here we choose $\zeta(\lambda)=\lambda^2(\lambda-1)^2$ and  $\alpha=0.05$. We also impose the pure Dirichlet boundary condition. 

We implement
the proposed method with    $\nu=1$ and $\nu=10^{-5}$, respectively. Numerical results are displayed in Table \ref{3D table}, which shows the optimal convergence orders and the pressure-robustness of the propsoed method in 3D.

 
\begin{table}[!htbp]
\centering
\caption{Example \ref{ex:3D}: Numerical results of the proposed method with $\nu=1$ and $\nu=10^{-5}$.}
\label{3D table}
\begin{tabular}{@{}lcccccc@{}}
\toprule
$h$ & $\left\|{\bm u}-{\bm u}_0\right\|$ & order & $|{\bm u}-{\bm u}_0|_1$ & order & $\left\|p-p_h\right\|$ & order \\
\midrule
\multicolumn{7}{c}{$\nu=1$} \\
\midrule
1/8  & $4.427\mathrm{e}-06$ &      & $8.787\mathrm{e}-04$ &      & $2.439\mathrm{e}-03$ &      \\
1/12 & $2.056\mathrm{e}-06$ & 1.89 & $5.989\mathrm{e}-05$ & 0.94 & $1.631\mathrm{e}-03$ & 0.99 \\
1/16 & $1.174\mathrm{e}-06$ & 1.94 & $4.526\mathrm{e}-05$ & 0.97 & $1.225\mathrm{e}-03$ & 0.99 \\
1/20 & $7.570\mathrm{e}-07$ & 1.96 & $3.633\mathrm{e}-05$ & 0.98 & $9.808\mathrm{e}-04$ & 0.99 \\
\midrule
\multicolumn{7}{c}{$\nu=10^{-5}$} \\
\midrule
1/8  & $4.422\mathrm{e}-06$ &      & $8.789\mathrm{e}-04$ &      & $2.439\mathrm{e}-03$ &      \\
1/12 & $2.053\mathrm{e}-06$ & 1.89 & $5.990\mathrm{e}-05$ & 0.94 & $1.631\mathrm{e}-03$ & 0.99 \\
1/16 & $1.172\mathrm{e}-06$ & 1.94 & $4.526\mathrm{e}-05$ & 0.97 & $1.225\mathrm{e}-03$ & 0.99 \\
1/20 & $7.557\mathrm{e}-07$ & 1.96 & $3.633\mathrm{e}-05$ & 0.98 & $9.807\mathrm{e}-04$ & 0.99 \\
\bottomrule
\end{tabular}
\end{table}

  

\end{example}

\begin{example}[No flow test]\label{No flow}\rm
To further verify the pressure-robustness  of the proposed method, we consider a no flow test. We choose 
the velocity and pressure as
$$
\bm{u}=\left(\begin{array}{l}
0 \\
0
\end{array}\right), \quad p=-\frac{\mathrm{Ra}}{2} y^2+\operatorname{Ra} y-\frac{\mathrm{Ra}}{3}
$$
where $\mathrm{Ra}=1000$.

We choose  $\nu=1$. In this test, the true velocity is 0. 
 Figure \ref{fig1} displays the two components of the numerical velocity, which  achieve machine precision.   This demonstrates that even in the presence of large pressure, there is no pollution on the numerical velocity, indicating the pressure-robustness of the
proposed method.
 
\begin{figure}[htbp]
\centering
 \subfigure[$u_1$] {
\includegraphics[width=0.48\textwidth]{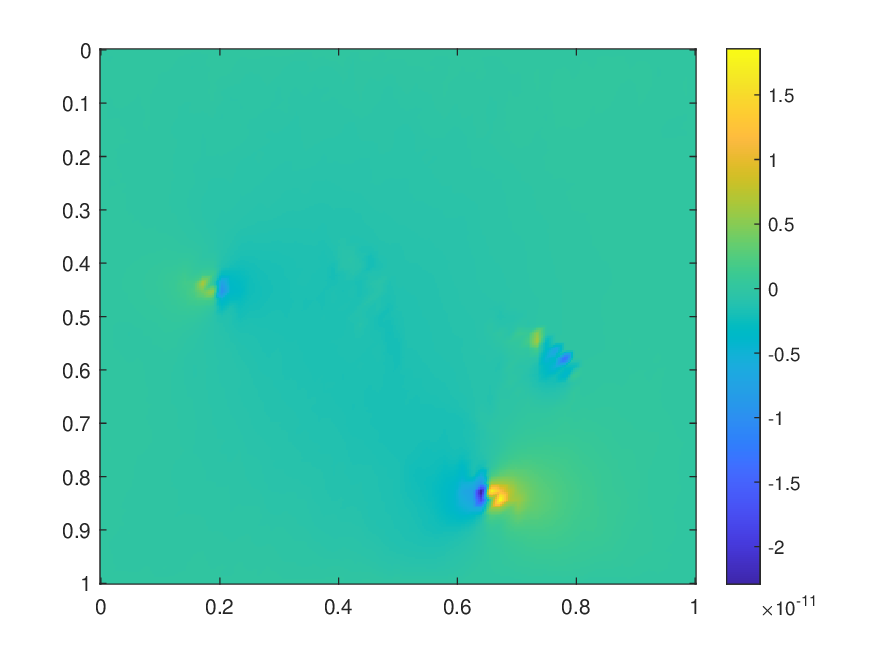}}
\subfigure[$u_2$]  {
\centering
\includegraphics[width=0.48\textwidth]{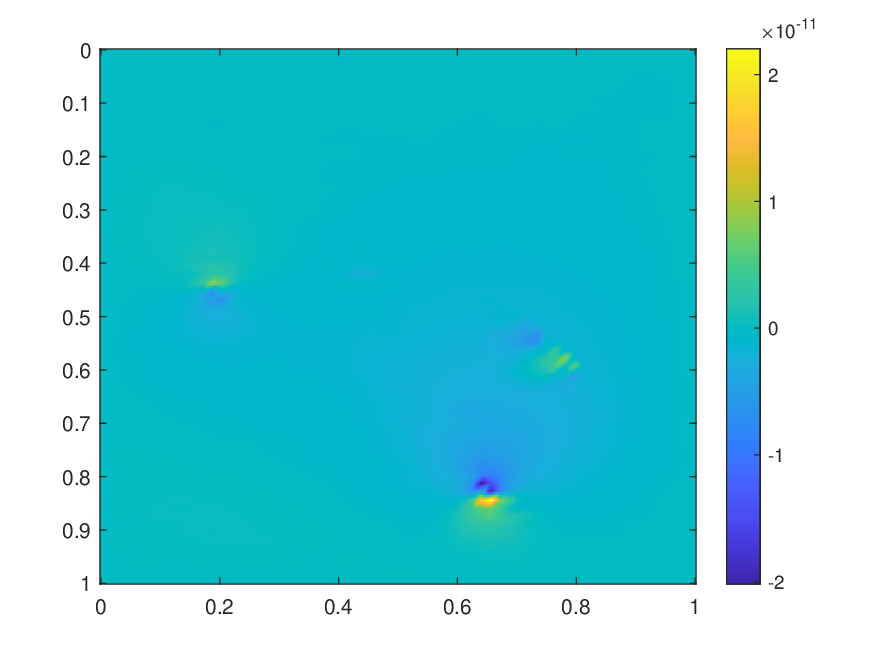}
 }\caption{Example \ref{No flow}: Plots of the numerical velocity.}
 \label{fig1}
\end{figure}
\end{example}



\begin{example}[Two-dimensional Lid-driven Cavity Flow]\label{Lid-driven}\rm
In this numerical example, we consider a lid-driven cavity flow over the domain $\Omega = (0,1)^2$ to verify that the computed velocity is independent of the irrotational component of the body force ($\nabla \times \bm{f} = 0$).
To this end, we set
\begin{itemize}
    \item the Dirichlet boundary condition:
    $$
    \left.\bm u\right|_{\partial \Omega}=\bm u_D= \begin{cases}(1,0)^{\mathrm T}, & \text { if } y=1, \\ (0,0)^{\mathrm T}, & \text { otherwise}, \end{cases}
    $$ 
    \item two distinct body forces: $ \bm{f}_1 = 0,\ \bm{f}_2 = \frac{10^6}{3} \nabla \left(x^3 + y^3\right)$ (see Figure \ref{fig_Ld}(a)).
\end{itemize}

First, setting $\nu = 1$, we solve the NS equations on a mesh with $h=1/100$ using the two different body forces $\bm f = \bm f_1$ and $\bm f = \bm f_2$.  Figures~\ref{fig_Ld}(b) and~\ref{fig_Ld}(c) display the differences in the first and second components of the resulting velocity fields, respectively, showing negligible discrepancies and thereby confirming the theoretical result stated in Theorem~\ref{th:error}.

{
Furthermore, with Reynolds number defined as $\operatorname{Re} = \frac{1}{\nu}$, we present streamline contours for the lid-driven cavity flow with body force $\bm f=\bm f_1$ for $\operatorname{Re}=5{,}000$, $15{,}000$, and $22{,}000$, in Figure \ref{figlid}. The computations remain stable across these Reynolds numbers and clearly resolve both the primary vortex and the corner vortices. 
In \cite{erturk2005numerical}, stability is reported up to \(\operatorname{Re}=21{,}000\) with \(h=1/601\), while a more recent work \cite{yang2022analysis} shows streamline contours at \(\operatorname{Re}=20{,}000\) using high-order elements on a boundary-refined mesh. By contrast, our low-order scheme achieves stable solutions up to \(\operatorname{Re}=22{,}000\) on a far coarser uniform mesh $h=1/250$ than \cite{erturk2005numerical}, evidencing the robustness and efficiency of our approach.

}
\begin{figure}[htbp]
\centering
\subfigure[] {  \includegraphics[width=0.31\textwidth]{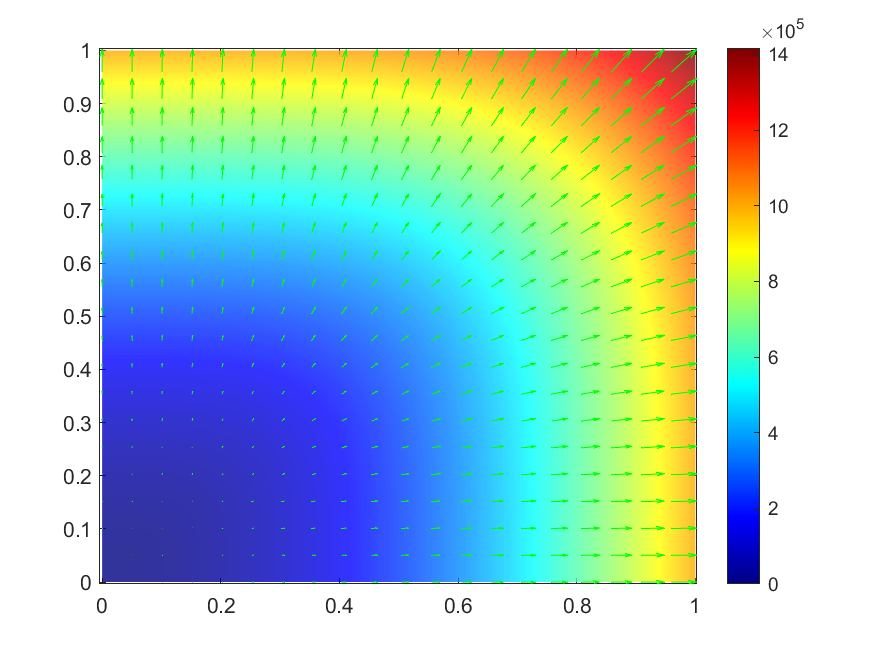}
}
\subfigure[] {
\includegraphics[width=0.31\textwidth]{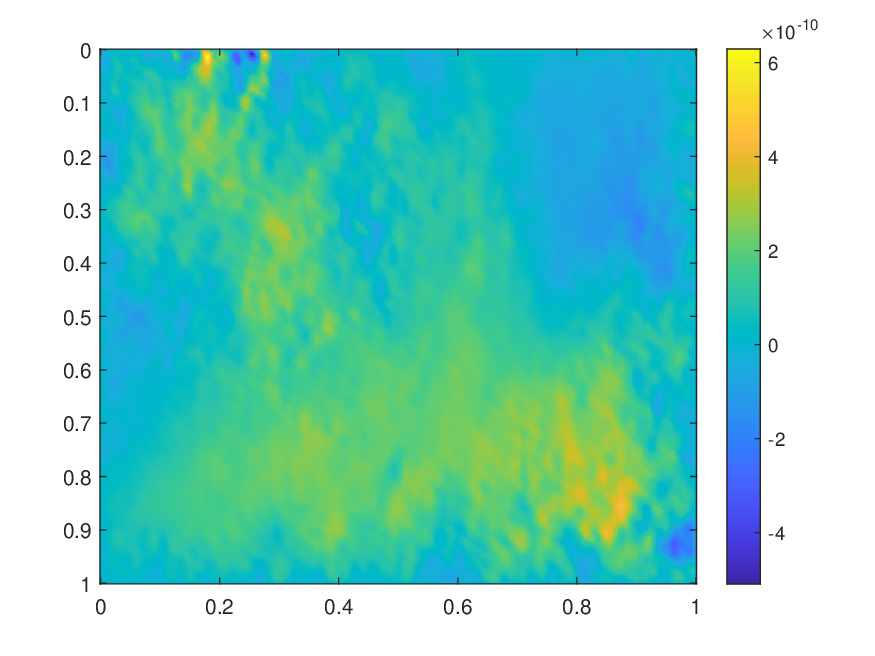}
}
\subfigure[] { \includegraphics[width=0.31\textwidth]{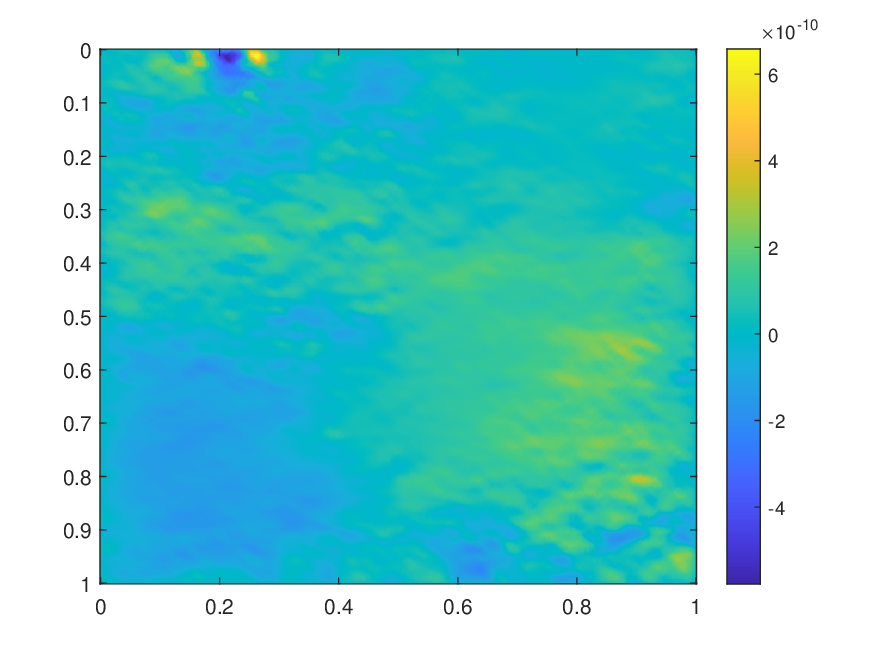}
}
\caption{Example \ref{Lid-driven}: (a) body force $\bm f_2$; (b) difference in the first velocity component computed with $\bm f_1$ and $\bm f_2$; (c) difference in the second velocity component computed with $\bm f_1$ and $\bm f_2$.}
\label{fig_Ld}
\end{figure}

\begin{figure}[htbp]
    \centering
\includegraphics[width=0.65\textwidth]{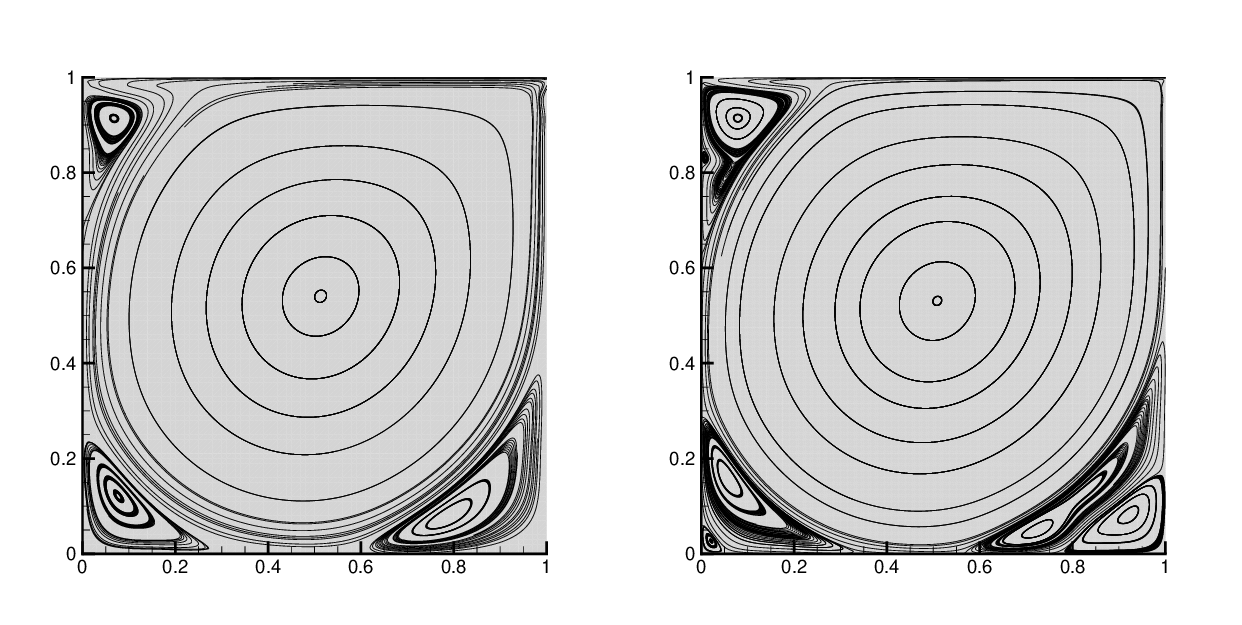}
\includegraphics[width=0.325\textwidth]{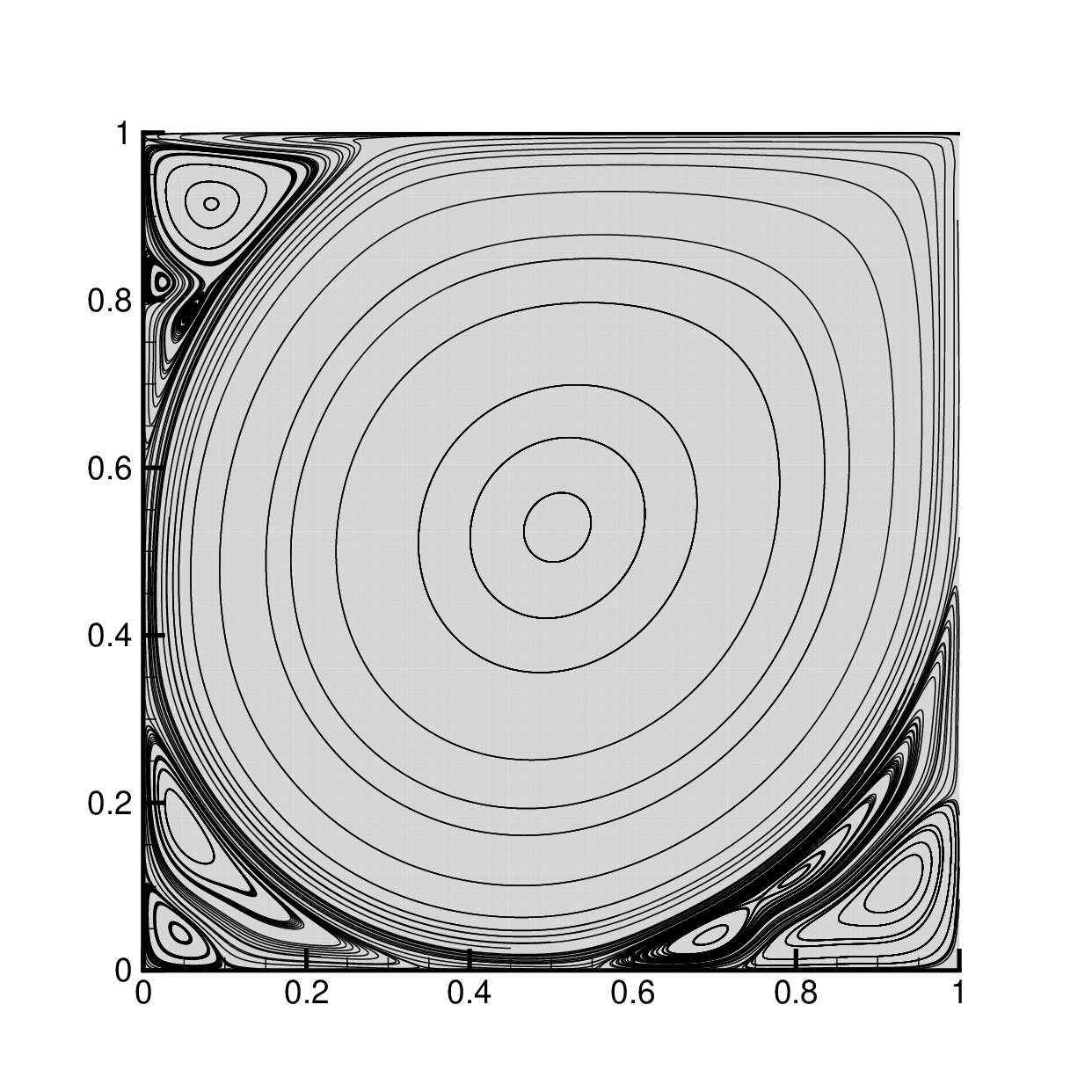}
\caption{Example \ref{Lid-driven}: Streamlines of lid-driven cavity flow at Re = 5000 (left) with $h = 1/100$,  Re = 15000 (middle) with $h = 1/200$, and Re = 22000 (right) with $h = 1/250$.}
\label{figlid}
\end{figure}


\end{example}

{\begin{example}[Backward facing step test]\label{BFS_example}\rm
In this example, we consider the backward facing step problem on an
L-shaped domain $\Omega = (-4, 20)\times(0, 2)\backslash[-4, 0] \times [0, 1]$ with $\bm f=0$. The boundary conditions are specified as follows:
\begin{itemize}
\item
two different inflow profiles on $x = -4$:
\begin{itemize}
\item parabolic inflow: $\bm{u} = (6(y - 1)(2 - y), 0)^{\mathrm{T}}$,
\item constant inflow: $\bm{u} = (1, 0)^{\mathrm{T}}$, 
\end{itemize}
\item
the outflow condition on $x = 20$:
$\left(\nu \nabla \bm{u} - p^{\text{kin}} \mathbb{I} \right) \bm{n} = \bm{0}$,
\item
homogeneous Dirichlet boundary conditions on the remaining boundaries:
$\bm{u} = \bm{0}$.
\end{itemize}

The corresponding Newton’s iteration algorithm of the method presented in Remark \ref{mixboundary} reads
\begin{align*}
a \left({\bm u}_h^{n+1},\bm v\right)&+ c \left({\bm u}_h^{n+1},{\bm u}_h^{n},\bm v\right)+c\left({\bm u}_h^{n},{\bm u}_h^{n+1},\bm v\right)-b \left(\bm v, p_h^{n+1}\right) +d\left({\bm u}_h^{n+1},{\bm u}_h^{n},\bm v\right)\\+&d\left({\bm u}_h^{n},{\bm u}_h^{n+1},\bm v\right)=\left({\bm f}, \mathcal{R} \bm v\right)+ c\left({\bm u}_h^{n},{\bm u}_h^{n},\bm v\right)+d\left({\bm u}_h^{n},{\bm u}_h^{n},\bm v\right),\quad \forall \bm v \in \bm V_h^{0,D}.
\end{align*}

Figures \ref{figbj1} and \ref{figbj2} depict 
the streamlines at \(\operatorname{Re} = \frac{1}{\nu}=100\), 500, and 1000, corresponding to the parabolic and constant inlet profiles, respectively. The flow topology exhibits similar patterns under both inlet conditions, with a primary recirculation zone forming downstream of the step and elongating as Re increases. A closer comparison reveals that, as shown in Figures  \ref{figbj1} and \ref{figbj2}, the reattachment length of the primary recirculation zone is significantly longer for the parabolic inlet condition than for the constant inlet condition. 
Comparison with the results in \cite{abide20052d} indicates that our computed streamlines for \(\operatorname{Re} = 100 \text{ and }\operatorname{Re} =500 \) are in good agreement with the literature, while the streamlines for \(\operatorname{Re} =1000 \), which are also presented here, provide an extension beyond the cases considered in \cite{abide20052d}.




\begin{figure}[htbp]
    \centering
\includegraphics[width=1\textwidth]{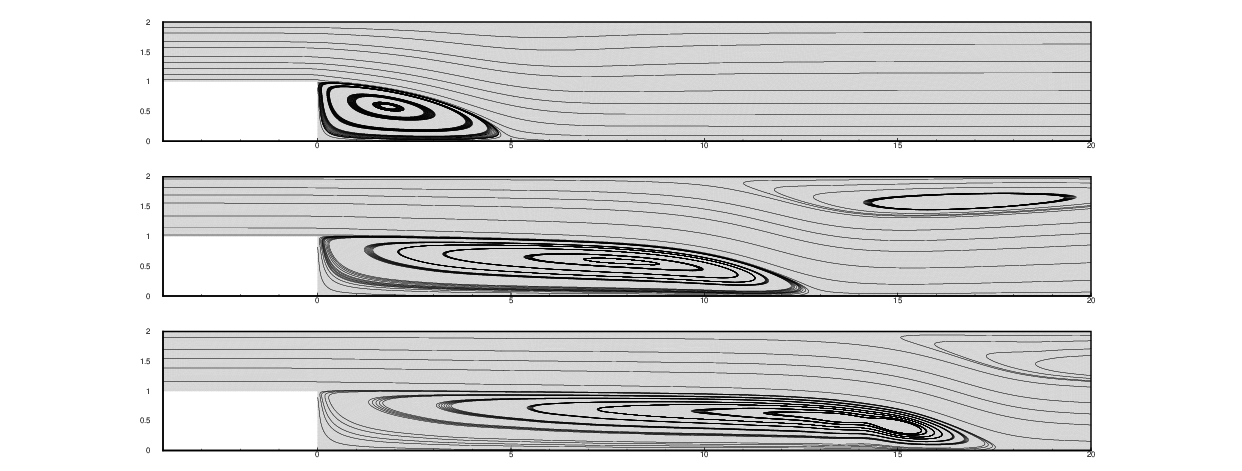}
\caption{Example \ref{BFS_example}: Streamlines of  backward facing step flows with parabolic inlet at Re = 100, 500, and 1000 from top to bottom. }
\label{figbj1}
\end{figure}

\begin{figure}[htbp]
    \centering
\includegraphics[width=1\textwidth]{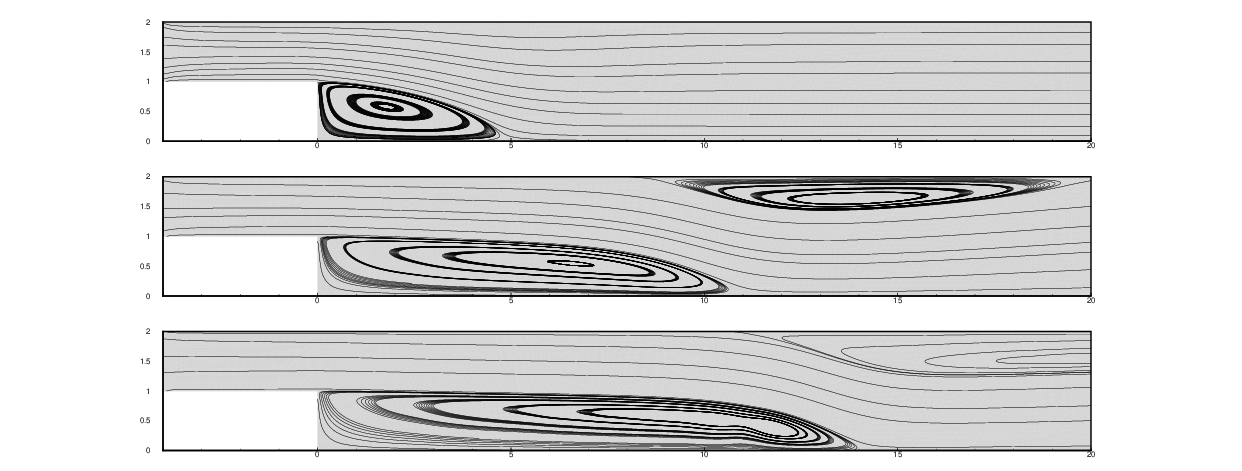}
\caption{Example \ref{BFS_example}: Streamlines of backward facing step flows with constant inlet at Re = 100, 500, and 1000 from top to bottom. }
\label{figbj2}
\end{figure}


\end{example}
}

\begin{example}[Laminar flow around a cylinder]\label{cylinder} \rm
In this example, we consider the famous Laminar flow past a cylinder \cite{schafer1996benchmark} to investigate the appearance and evolution of the symmetric eddies behind the cylinder as the Reynolds number increases.
Let the domain $ \Omega=(0,2.2) \times(0,0.41) \backslash \overline{B_{0.05}(0.2,0.2)}$
 where $B_r\left(x, y\right)$ denotes a ball with center $\left(x, y\right)$ and radius $r$. Specifically, 
we set
 \begin{itemize}
     \item the body force $\bm{f}=\mathbf{0}$,
      \item the
no-slip boundary conditions at the walls $y=0$, $y=0.41$, $\partial B_{0.05}(0.2,0.2)$: $\bm{u}=\mathbf{0}$,
       \item 
the inlet boundary condition on $x=0$: $\bm{u}=(1,0)^{\mathrm{T}}$,
        \item the outlet boundary condition on $x=2.2$:
        $\left(\nu \nabla \bm{u}-p^{\text {kin}} \mathbb I\right) \bm{n}=\mathbf{0}$. 
        \end{itemize}


As in \cite{yang2022analysis}, we calculate the Reynolds number using the formula 
\[
\operatorname{Re} = \frac{2 \times 0.05 \times 1}{\nu} = \frac{1}{10 \nu},
\]
and set \(\operatorname{Re} = 5,\ 10,\ 40,\text{ and }100\).

The contours of velocity magnitude $|\bm u_0|$ and kinematic pressure $p^{\mathrm{kin}}_h = p_h - \frac{1}{2} \mathcal Q_h(|\bm u_0|^2)$ for various Reynolds numbers are illustrated in  Figures \ref{fig:contourv} and \ref{fig:contourp}, respectively. These figures clearly demonstrate robust symmetry, and as the Reynolds number increases, regions of high velocity magnitude extend towards the outlet, consistent with findings in \cite{franca2001two,yang2022analysis}.

Furthermore, Figure \ref{fig:streamlice} presents a detailed depiction of the streamlines at different Reynolds numbers around the circulation area. Symmetric eddies appear behind the cylinder starting at $\operatorname{Re}=10$, aligning well with the observations in \cite{van1975perturbation}. With increasing Reynolds numbers, the two eddies enlarge and shift towards the outlet while maintaining symmetry, corroborating the results in \cite{yang2022analysis}.

\begin{figure}[htbp]
    \centering
\includegraphics[width=1\textwidth]
{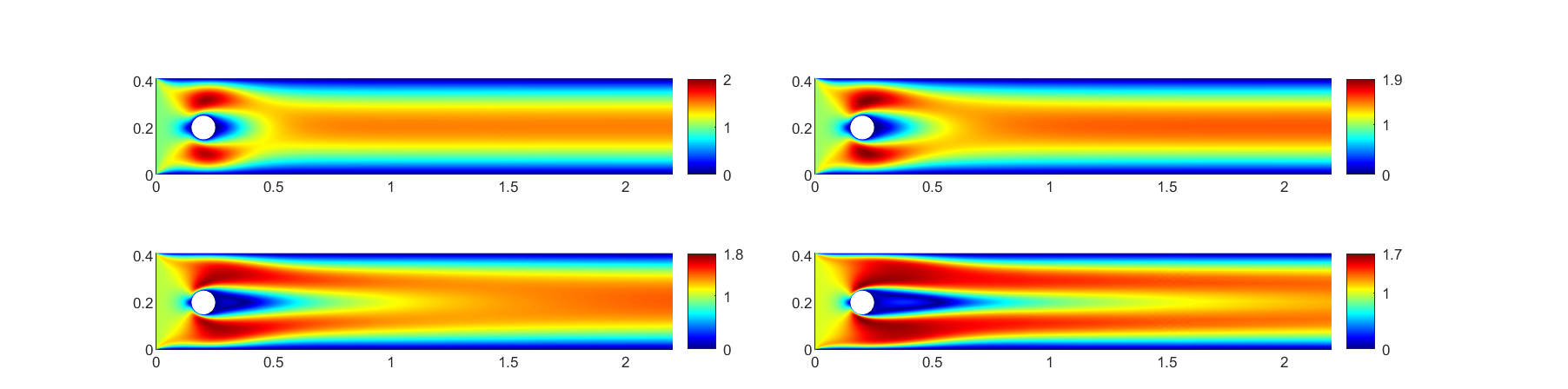}
    \caption{Example \ref{cylinder}: Contours of velocity magnitude for Re = 5, 10, 40, and 100 from top-left to bottom-right.}
    \label{fig:contourv}
\end{figure}

\begin{figure}[htbp]
    \centering
\includegraphics[width=1\textwidth]{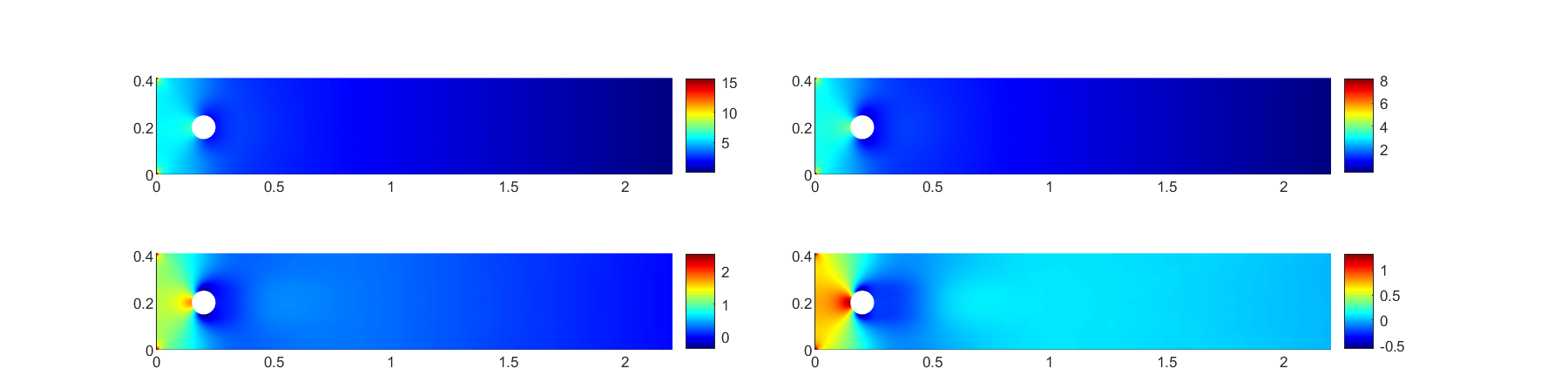}
    \caption{Example \ref{cylinder}: Contours of kinematic pressure for Re = 5, 10, 40, and 100 from top-left to bottom-right.}
    \label{fig:contourp}
\end{figure}

\begin{figure}[htbp]
    \centering
\includegraphics[width=1\textwidth]{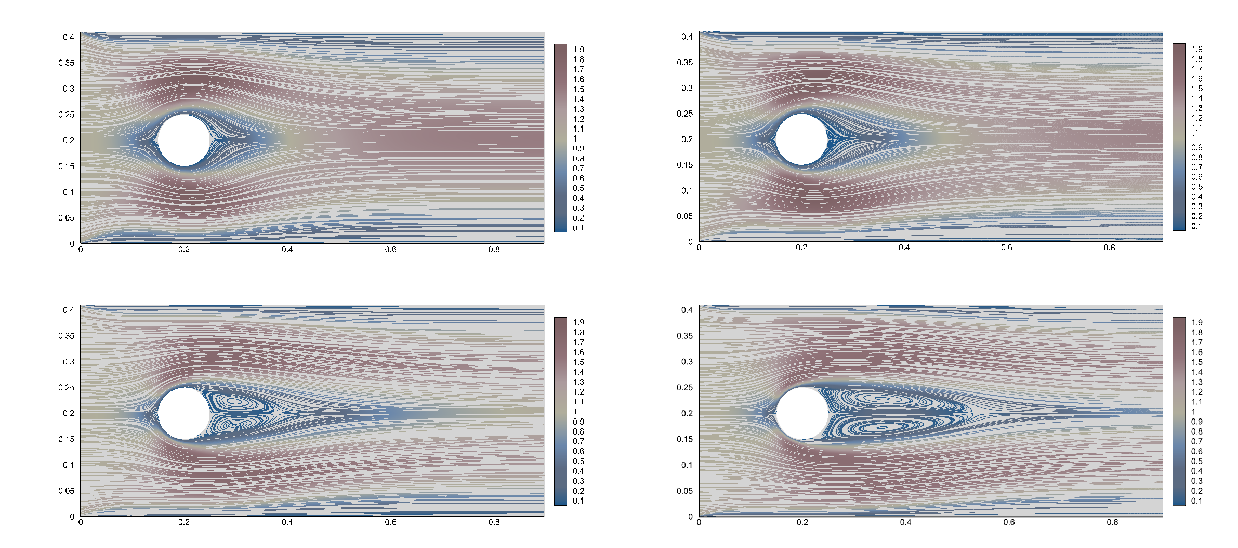}
    \caption{Example \ref{cylinder}: Streamlines on colored velocity magnitude distribution in the circulation area for Re = 5, 10, 40, and 100 from top-left to bottom-right.}
    \label{fig:streamlice}
\end{figure}

\end{example}

\section{Conclusion}\label{conclusion}
In this paper, we developed a PR\&PF-EG method for the steady incompressible NS in rotational form. 
We rigorously establish the well-posedness of the proposed method and derive pressure-robust error estimates. The numerical experiments validate the theoretical results and demonstrate the robustness and efficiency of our method, especially at high Reynolds numbers.

The method is based on an EG space that enriches the first-order CG space with piecewise constants on edges or faces. It is particularly well-suited for the rotational formulation, as the curl of an EG function is entirely determined by its CG component and remains unaffected by the DG enrichment. 
Moreover, \(\nabla \times\bm u_h = \nabla \times\bm u_0\) falls within the space of \(\operatorname{RT}_0\), aligning with the image of the velocity reconstruction operator. Consequently, it is feasible to choose a $v_b$ in test function \(\bm v_h=\{\bm v_0,v_b\}\) such that \(\mathcal R\bm v_h = \nabla \times\bm u_0\) to cancel the nonlinear term $(\nabla\times\bm u_0\times\mathcal R\bm u_h,\mathcal R\bm v_h)$. This characteristic is crucial for the development of a helicity-conserving numerical scheme \cite{rebholz2007energy}. Future work will explore extending this method to a helicity-conserving EG method for unsteady NS equations.


\section*{Declarations}
\bmhead{Funding}

This work was partially  supported by the National Natural Science Foundation of
China (No. 12201020).

\bmhead{Conflict of Interest}

The authors declare that they have no conflict of interests.




\bmhead{Code availability}
The code used in this work will be made available upon reasonable request.

\bmhead{Author contribution} 
Conceptualization: Qian Zhang; Methodology: Shuai Su, Qian Zhang; Investigation: Shuai Su, Xiurong Yan, Qian Zhang;
Software: Xiurong Yan, Qian Zhang; Writing - original draft preparation: Shuai Su, Xiurong Yan; Writing - review and editing: Shuai Su, Qian Zhang; Funding acquisition: Shuai Su; Project administration: Shuai Su, Qian Zhang; Supervision: Shuai Su, Qian Zhang.

\bibliography{sn-bibliography}

\end{document}